\newcommand{\Fisher}{Fisher:1948-short}

\documentclass[10pt]{article}
\usepackage{amsmath,amsthm,amsfonts,amssymb,latexsym,graphicx,subfig,mathtools,multirow}
\usepackage{algorithm}
\usepackage{algpseudocode}
\usepackage[pdfpagemode=UseNone,pdfstartview=FitH]{hyperref}
\newcommand{\Extra}[1]{}

\usepackage{comment}

\allowdisplaybreaks[4]

\newcommand{\R}{\mathbb{R}}
\DeclareMathOperator{\Prob}{\mathbb{P}}
\DeclareMathOperator{\Expect}{\mathbb{E}}

\newcommand{\E}{\Expect}               
\newcommand{\q}{q}
\newcommand{\pretend}[2]{\smash{\mathrlap{#1}}\phantom{#2}}
\newcommand{\qwc}{\pretend{\underline{q}}{q}}

\newcommand{\w}{\mathbf{w}}

\newcommand{\BA}{\mathrm{BA}}
\newcommand{\BG}{\mathrm{BG}}

\DeclareMathOperator{\III}{\boldsymbol{1}}

\newcommand{\st}{\,\middle|\,}
\newcommand{\givn}{\mathrel{|}}
\newcommand{\dd}{\mathrm{d}}
\renewcommand{\complement}{\mathsf{c}}

\newtheorem{InnerCustomTheorem}{Theorem}
\newenvironment{CustomTheorem}[1]
  {\renewcommand\theInnerCustomTheorem{#1}\InnerCustomTheorem}
  {\endInnerCustomTheorem}
\newtheorem{InnerCustomProposition}{Proposition}
\newenvironment{CustomProposition}[1]
  {\renewcommand\theInnerCustomProposition{#1}\InnerCustomProposition}
  {\endInnerCustomProposition}
\newtheorem{theorem}{Theorem}
\newtheorem{lemma}[theorem]{Lemma}
\newtheorem{proposition}[theorem]{Proposition}

\theoremstyle{definition}
\newtheorem{remark}{Remark}
\usepackage[normalem]{ulem}

\title{Combining p-values via averaging}
\author{Vladimir Vovk\thanks%
    {Department of Computer Science,
    Royal Holloway, University of London,
    Egham, Surrey, UK.
    E-mail: \href{mailto:v.vovk@rhul.ac.uk}{v.vovk@rhul.ac.uk}.}
  \and Ruodu Wang\thanks%
    {Department of Statistics and Actuarial Science,
    University of Waterloo,
    Waterloo, Ontario, Canada.
    E-mail: \href{mailto:wang@uwaterloo.ca}{wang@uwaterloo.ca}.}}

\begin{document}
\maketitle

\begin{abstract}
  This paper proposes general methods for the problem of multiple testing of a single hypothesis,
  with a standard goal of combining a number of p-values without making any assumptions
  about their dependence structure.
  An old result by R\"uschendorf and, independently, Meng implies that the p-values can be combined
  by scaling up their arithmetic mean by a factor of 2
  (and no smaller factor is sufficient in general).
  A similar result about the geometric mean (Mattner) replaces 2 by $e$.
  Based on more recent developments in mathematical finance, specifically, robust risk aggregation techniques,
  we extend these results to generalized means;
  in particular, we show that $K$ p-values can be combined
  by scaling up their harmonic mean by a factor of $\ln K$
  (asymptotically as $K\to\infty$). 
  This leads to a generalized version of the Bonferroni--Holm procedure. 
  We also explore methods using weighted averages of p-values. 
  Finally, we  discuss the efficiency of various methods of combining p-values
  and how to choose a suitable method in light of data and prior information.
\end{abstract}

\section{Introduction}
\label{sec:introduction}

Suppose we are testing the same hypothesis using $K\ge2$ different statistical tests
and obtain p-values $p_1,\dots,p_K$.
How can we combine them into a single p-value?

Some early answers to this question include 
Tippett's \cite{Tippett:1931}, Pearson's \cite{Pearson:1933}, and Fisher's \cite{\Fisher}.
These solutions assume that the p-values are independent
(in practice, obtained from independent test statistics),
whereas we would like to avoid any assumptions besides all $p_1,\dots,p_K$
being bona fide p-values.
Fisher's method has been extended to dependent p-values in, e.g.,
\cite{Brown:1975,Kost/McDermott:2002},
but the combined p-values obtained in those papers are approximate;
in this paper we are interested only in precise or conservative p-values.
A summary of combination methods for p-values can be found in, e.g.,
Oosterhoff \cite{oosterhoff1976combination} and Cousins \cite{cousins2007annotated};
see also the more specialized review Owen \cite{owen2009karl}.

Without assuming any particular dependence structure among p-values, 
the simplest way of combining them is the Bonferroni method
\begin{equation}\label{eq:Bonferroni}
  F(p_1,\dots,p_K)
  :=
  K\min(p_1,\dots,p_K)
\end{equation}
(when $F(p_1,\dots,p_K)$ exceeds $1$ it can be replaced by $1$,
but we usually ignore this trivial step).
Albeit $F(p_1,\dots,p_K)$ is a p-value
(see Section \ref{sec:merging} for a precise definition of a \emph{p-value}),
it has been argued that in some cases it is overly conservative.
R\"uger \cite{Ruger:1978} extends the Bonferroni method by showing that,
for any fixed $k\in\{1,\dots,K\}$,
\begin{equation}\label{eq:Ruger}
  F(p_1,\dots,p_K)
  :=
  \frac{K}{k}p_{(k)}
\end{equation}
is a p-value,
where $p_{(k)}$ is the $k$th smallest p-value among $p_1,\dots,p_K$;
see \cite{Morgenstern:1980} for a simpler exposition.
Hommel~\cite{Hommel:1983} develops this by showing that
\begin{equation}\label{eq:Hommel}
  F(p_1,\dots,p_K)
  :=
  \left(
    1 + \frac12 + \dots + \frac1K
  \right)
  \min_{k=1,\dots,K}
  \frac{K}{k}p_{(k)}
\end{equation}
is also a p-value.
In the case of independent $p_1, \dots, p_n$,
Simes \cite{Simes:1986} improves \eqref{eq:Hommel} by removing the first factor
on the right-hand side of \eqref{eq:Hommel}.

A natural way to combine $K$ p-values is to average them, by using $\bar p:=(p_1+\dots+p_K)/K$
(at least when the tests have similar power).
Unfortunately, $\bar p$ is not necessarily a p-value. 
An old result by R\"uschendorf \cite[Theorem~1]{Ruschendorf:1982}
shows that $2\bar p$ is a p-value;
moreover, the factor of 2 cannot be improved in general.
In the statistical literature this result was rediscovered
by Meng \cite[Theorem~1]{Meng:1994}.

In this paper (see Section~\ref{sec:symmetric}) we move on to a general notion of the mean
as axiomatized by Andrei Kolmogorov \cite{Kolmogorov:1930mean}
and adapt various results of \emph{robust risk aggregation}
\cite{Embrechts/Puccetti:2006,Bernard/etal:2014,Embrechts/etal:2015,Wang/Wang:2016,Jakobsons/etal:2016}
to combining p-values by averaging them in Kolmogorov's wider sense.
In particular, to obtain a p-value from given p-values $p_1,\dots,p_K$,
it is sufficient to multiply their geometric mean by $e$
(as noticed by Mattner \cite{Mattner:2012-journal})
and to multiply their harmonic mean by $e\ln K$ (for $K>2$).
More generally, we consider the mean $M_{r,K}(p_1,\dots,p_K)$
(sometimes referred to as \emph{generalized mean})
defined by $((p_1^r+\dots+p_K^r)/K)^{1/r}$ for $r\in[-\infty,\infty]$,
and derive values of  $a_{r,K}$ making $a_{r,K}M_{r,K}$ a merging function.
In particular, our results cover the Bonferroni method~\eqref{eq:Bonferroni},
which corresponds to $M_{-\infty,K}(p_1,\dots,p_K)=K\min(p_1,\dots,p_K)$
(see, e.g., \cite[(2.3.1)]{Hardy/etal:1952}).

Median is also sometimes regarded as a kind of average.
R\"uger's~\eqref{eq:Ruger}, applied to $k:=\lceil K/2\rceil$, says that p-values
can be combined by scaling up their median by a factor of 2
(exactly for even $K$ and approximately for large odd $K$).
Therefore, we have the same factor of 2 as in R\"uschendorf's \cite{Ruschendorf:1982} result.
(Taking $k=\lfloor(K+1)/2\rfloor=\lceil K/2\rceil$ is suggested
in \cite[Section~1.1]{Mattner:2011-journal}.)
More generally, the $\alpha$ quantile $p_{(\lceil\alpha K\rceil)}$ becomes a p-value if multiplied by $1/\alpha$.

Section~\ref{sec:weighted}--\ref{sec:choice} barely scratch the surface of three important topics:
weighted averaging of p-values,
the efficiency of our methods of combining p-values under various alternative distributions,
and the choice of a suitable value of the parameter $r$ for generalized mean.
Section~\ref{sec:conclusion} concludes the paper listing some directions of further research.

As this paper targets a statistical audience,
in the main text we shall omit proofs and techniques based on results
from the literature on robust risk aggregation;
the details can be found in Appendix~\ref{app:proofs-validity}.

It is often possible to automatically transform
results about multiple testing of a single hypothesis
into results about testing multiple hypotheses;
the standard procedures are Marcus et al.'s \cite{Marcus/etal:1976} closed testing procedure
and its modification by Hommel~\cite{Hommel:1986}.
In particular, when applied to the Bonferroni method
the closed testing procedure gives the well-known procedure due to Holm \cite{Holm:1979},
which we will refer to as the Bonferroni--Holm procedure;
see, e.g., \cite{Hommel:1986,Hommel:1988} for its further applications.
In Appendix~\ref{sec:multiple} we briefly discuss a similar application
to one of the methods given in Section~\ref{sec:symmetric}.

\subsection*{Some terminology}

A function $F:[0,1]\to[0,\infty)$ is \emph{increasing} (resp.\ \emph{decreasing})
if $F(x_1)\le F(x_2)$ (resp.\ $F(x_1)\ge F(x_2)$) whenever $x_1\le x_2$.
A function $F:[0,1]^K\to[0,\infty)$ is \emph{increasing} (resp.\ \emph{decreasing})
if it is increasing (resp.\ decreasing) in each of its arguments.
A function is \emph{strictly increasing} or \emph{strictly decreasing}
when these conditions hold with strict inequalities.

\section{Merging functions}
\label{sec:merging}

A \emph{p-variable} is a random variable $P$ that satisfies
\begin{equation}\label{eq:p}
  \Prob(P\le\epsilon)
  \le
  \epsilon,
  \quad
  \forall\epsilon\in(0,1).
\end{equation}
The values taken by a p-variable are \emph{p-values}
(allowed to be conservative).
In Section~\ref{sec:introduction} the expression ``p-value''
was loosely used to refer to p-variables as well.
A \emph{merging function} is an increasing Borel function $F:[0,1]^K\to[0,\infty)$
such that $F(U_1,\dots,U_K)$ is a p-variable for any choices of  random variables $U_1,\dots,U_K$
(on the same probability space, which can be arbitrary)
distributed uniformly on $[0,1]$.
Without loss of generality we can assume that $U_1,\dots,U_K$ are defined on the same atomless probability space,
which is fixed throughout the paper
(cf.\ \cite[Proposition~A.27]{Follmer/Schied:2011}).
Let $\mathcal{U}$ be the set of all uniformly distributed random variables
(on our probability space).
Using the notation $\mathcal{U}$,
an increasing Borel function $F:[0,1]^K\to[0,\infty)$ is a merging function if, for each $\epsilon\in(0,1)$,
\begin{equation}\label{eq:merging}
  \sup
  \left\{
    \Prob(F(U_1,\dots,U_K)\le\epsilon)
    \st
    U_1,\dots,U_K\in\mathcal{U}
  \right\}
  \le
  \epsilon.
\end{equation}
We say that 
a merging function $F$ is \emph{precise} if, for each $\epsilon\in(0,1)$,
\begin{equation}\label{eq:precise-merging}
  \sup
  \left\{
    \Prob(F(U_1,\dots,U_K)\le\epsilon)
    \st
    U_1,\dots,U_K\in\mathcal{U}
  \right\}
  =
  \epsilon.
\end{equation}

\begin{remark}
  The requirement that a merging function be Borel
  does not follow automatically from the requirement that it be increasing:
  see the remark after Theorem~4.4 in \cite{Graham/Grimmett:2006}
  (Theorem~4.4 itself says that every increasing function on $[0,1]^K$ is Lebesgue measurable).
\end{remark}

It may be practically relevant to notice that, for any merging function $F$,
$F(P_1,\dots,P_K)$ is a p-variable whenever $P_1,\dots,P_K$ are p-variables
(on the same probability space).
Indeed, for each $k\in\{1,\dots,K\}$ we can define a uniformly distributed
(see, e.g., \cite[Proposition 2.1]{ruschendorf2009distributional})
random variable $U_k\le P_k$ by
\[
  U_k(\omega)
  :=
  \Prob(P_k<P_k(\omega))
  +
  \Theta(\omega)
  \Prob(P_k=P_k(\omega)),
  \enspace
  \omega\in\Omega,
\]
where $\Theta$ is a random variable distributed uniformly on $[0,1]$ and independent of $P_1,\dots,P_K$,
and $\Omega$ is the underlying probability space extended (if required)
to carry such a $\Theta$;
we then have
\[
  \Prob(F(P_1,\dots,P_K)\le\epsilon)
  \le
  \Prob(F(U_1,\dots,U_K)\le\epsilon)
  \le
  \epsilon,
  \quad
  \forall\epsilon\in(0,1).
\]
Therefore, combining p-values can be carried out in multiple layers
(although it may make the final combined p-value overly conservative);
we will discuss this further in Subsection~\ref{subsec:search}.

\section{Combining p-values by symmetric averaging}
\label{sec:symmetric}

In this section we present our methods of combining p-values via averaging. 
A general notion of averaging,
axiomatized by Kolmogorov \cite{Kolmogorov:1930mean}, is
\begin{equation}\label{eq:mean}
  M_{\phi,K}(p_1,\dots,p_K)
  :=
  \psi
  \left(
    \frac{\phi(p_1)+\dots+\phi(p_K)}{K}
  \right),
\end{equation}
where $\phi:[0,1]\to[-\infty,\infty]$ is a continuous strictly monotonic function and $\psi$ is its inverse
(with the domain $\phi([0,1])$).
For example, arithmetic mean corresponds to the identity function $\phi(p)=p$,
geometric mean corresponds to $\phi(p)=\ln p$, and harmonic mean corresponds to $\phi(p)=1/p$.

The problem of finding precise p-value merging functions based on the averaging method in \eqref{eq:mean}
is closely connected to \emph{robust risk aggregation},
an active topic in  mathematical finance
(see Remark~\ref{rem:qrm} below).
The origin of this field lies in a problem posed by Kolmogorov (see, e.g., \cite{Makarov:1981})
about bounds on the distribution function of the sum of random variables.
To establish the main results in this paper, we will use many recent results in robust risk aggregation, in particular,
many results in \cite{Embrechts/Puccetti:2006,Wang/etal:2013,Bernard/etal:2014,Embrechts/etal:2015,Wang/Wang:2016,Jakobsons/etal:2016}.

Below, we first briefly discuss the connection between the two problems in Section \ref{sec:31},
and then present our main results in Section \ref{sec:32}.

\subsection{Quantiles and robust risk aggregation}\label{sec:31}

We start from a simple result (Lemma~\ref{lem:RW-reduction} below)
that translates probability statements about merging functions into corresponding statements about quantiles.
This result will allow us to freely use some recent results in the literature on robust risk aggregation.
The proof of Lemma~\ref{lem:RW-reduction} is given on p.~\pageref{pf:lem1} in Appendix~\ref{app:proofs-validity}.

Define the left $\alpha$-quantile of a random variable $X$ for $\alpha\in(0,1]$,
\[
  \q_\alpha (X)
  :=
  \sup \{x \in \R: \Prob(X\le x)<\alpha\},
\]
and the right $\alpha$-quantile of $X$ for $\alpha\in[0,1)$,
\[
  \q^+_\alpha (X)
  :=
  \sup \{x \in \R: \Prob(X\le x)\le \alpha\}.
\]
Notice that $\q_1(X)$ is the essential supremum of $X$ and $\q_0^+(X)$ is the essential infimum of $X$.
For a function $F:[0,1]^K\to[0,\infty)$ and $\alpha \in (0,1)$,
write
\[
  \qwc_{\alpha}(F)
  :=
  \inf
  \left\{
    \q_\alpha(F(U_1,\dots,U_K))
    \st
    U_1,\dots,U_K\in\mathcal{U}
  \right\}.
\]

\begin{lemma}\label{lem:RW-reduction}
  For an increasing Borel function $F:[0,1]^K\to[0,\infty)$:
  \begin{itemize}
  \item[(a)]
    $F$ is a merging function if and only if
    $\qwc_\epsilon(F)\ge \epsilon$ for all $\epsilon\in(0,1)$;
  \item[(b)]
    $F$ is a precise merging function if and only if
    $\qwc_\epsilon(F)=\epsilon$ for all $\epsilon\in(0,1)$.
  \end{itemize}
\end{lemma}
 
\begin{remark}\label{rem:qrm}
  This remark discusses how the problem of combining p-values is related to robust risk aggregation
  in the field of mathematical finance.
  In quantitative risk management,
  the term \emph{robust risk aggregation} refers to evaluating
  the value of a \emph{risk measure} of an aggregation of risks $X_1,\dots,X_K$
  with specified marginal distributions and unspecified dependence structure.
  More specifically, if the risk measure is chosen as a quantile $q_\alpha$,
  known as \emph{Value-at-Risk} and very popular in finance,
  the quantities of interest are typically
  \[
    \overline q
    :=
    \sup
    \left\{
      \q_\alpha(X_1+\dots+X_n)
      \st
      X_1\sim F_1,\dots, X_n\sim F_n
    \right\}
  \]
  and
  \[
    \underline q
    :=
    \inf
    \left\{
      \q_\alpha(X_1+\dots+X_n)
      \st
      X_1\sim F_1,\dots, X_n\sim F_n
    \right\},
  \]
  where $F_1,\dots,F_n$ denote the prespecified marginal distributions of the risks.
  The motivation behind this problem is that, in practical applications of banking and insurance,
  the dependence structure among risks to aggregate is very difficult to accurately model,
  as compared with the corresponding marginal distributions.
  The interval $[\underline q, \overline q]$ thus represents all possible values of the aggregate risk measure given the marginal information.
  A more detailed introduction to this topic can be found in \cite[Section 8.4.4]{McNeil/etal:2015}
  and \cite[Chapter 4]{Ruschendorf:2013}.
  Via Lemma~\ref{lem:RW-reduction},
  the quantities $\overline q$ and $\underline q$ are obviously closely related to the problem of combining p-values.
  There are few explicit formulas for $\overline q$ and $\underline q$
  but fortunately some do exist in the literature, and they become useful in our study of merging functions.
\end{remark}

\subsection{Main results}\label{sec:32}

In this section, we present our main results on merging functions. The proofs of   these results  
along with a few auxiliary results,
can be found in Appendix~\ref{app:proofs-validity}.
We mostly focus on an important special case of \eqref{eq:mean}, namely,
\begin{equation}\label{eq:mean-r}
  M_{r,K}(p_1,\dots,p_K)
  :=
  \left(
    \frac{p_1^r+\dots+p_K^r}{K}
  \right)^{1/r},
\end{equation}
where $r\in\R\setminus\{0\}$ and the following standard conventions are used:
$0^c:=\infty$ for $c<0$, $0^c:=0$ for $c>0$,
$\infty+c:=\infty$ for $c\in\R\cup\{\infty\}$, and $\infty^c:=0$ for $c<0$.
The case $r=0$ (considered in \cite{Mattner:2012-journal})
is treated separately (as the limit as $r\to0$):
\begin{equation*}
  M_{0,K}(p_1,\dots,p_K)
  :=
  \exp
  \left(
    \frac{\ln p_1+\dots+\ln p_K}{K}
  \right)
  =
  \left(\prod_{k=1}^K p_k\right)^{1/K},
\end{equation*}
where, as usual, $\ln 0:=-\infty$, $-\infty+c:=-\infty$ for $c\in\R\cup\{-\infty\}$, and $\exp(-\infty):=0$.
It is also natural to set
\begin{align*}
  M_{\infty,K}(p_1,\dots,p_K)
  &:=
  \max(p_1,\dots,p_K),\\
  M_{-\infty,K}(p_1,\dots,p_K)
  &:=
  \min(p_1,\dots,p_K).
\end{align*}
The most important special cases of $M_{r,K}$ are perhaps
those corresponding to $r=-\infty$ (minimum), $r=-1$ (harmonic mean), $r=0$ (geometric mean),
$r=1$ (arithmetic mean), and $r=\infty$ (maximum);
the cases $r\in\{-1,0,1\}$ are known as Platonic means.

Our main aim is to identify merging functions of the form
\begin{equation}
  a_{r,K}
  M_{r,K}(p_1,\dots,p_K),~~ r\in [-\infty,\infty],~K=2,3,\dots,
  \label{eq:r1-3}
\end{equation}
where $a_{r,K}$ is a constant, hopefully making the merging function \eqref{eq:r1-3} precise.
In cases where an explicit formula of $a_{r,K}$ for \eqref{eq:r1-3} to be precise is not available, we obtain an asymptotically sharp bound  in explicit form.  
The  main results are summarized in Table~\ref{tab:main},
where a family $F_K$, $K=2,3,\dots$, of merging functions on $[0,1]^K$ is called \emph{asymptotically precise}
if, for any $a\in(0,1)$, the function $a F_K$ is not a merging function for a large enough $K$;
in other words, this family of merging functions cannot be improved by a constant multiplier.
It is well known \cite[Theorem~16]{Hardy/etal:1952} that 
\begin{equation} \label{eq:r1-1} 
 M_{r_1,K}\le M_{r_2,K}  \mbox{~on~}  [0,1]^K  \mbox{~if~}  r_1\le r_2 ;
\end{equation}
therefore, the constant $a_{r,K}$ making \eqref{eq:r1-3} a precise merging function should be generally decreasing in $r$.

\begin{table}[htbp]
\caption{The main results of Section~\ref{sec:symmetric}: examples of merging functions,
  all of them precise or asymptotically precise 
  (except for the case $r=-1$ where the asymptotic formula is not a merging function for finite $K$)}
\label{tab:main} 
\renewcommand{\arraystretch}{1.5}
\hspace*{-1cm}
\begin{tabular}{l|l|l|l|l}
  \hline 
  range of $r$ & \multicolumn{2}{l|}{merging function} & special case & claimed in\\
  \hline
  $r=\infty$ & $M_{r,K}$ & precise & maximum & \\\hline
  $r\in[K-1,\infty)$ & $K^{1/r}M_{r,K}$& precise& &  Prop.~\ref{prop:super-precise}\\\hline
  $r\in[\frac{1}{K-1},K-1]$ & $(r+1)^{1/r}M_{r,K}$& precise & arithmetic  &  Prop.~\ref{prop:precise}\\\hline
  $r\in(-1,\infty] $ & $(r+1)^{1/r}M_{r,K}$ & asymptotically precise && Prop.~\ref{prop:asymptotically-precise}\\\hline
  $r=0$ & $e M_{r,K}$ & asymptotically precise & geometric &   Prop.~\ref{prop:geometric}\\\hline
  \multirow{2}{4em}{$r=-1$}  & $e (\ln K)M_{r,K}$  &  not precise &   \multirow{2}{4em}{harmonic}  &  \multirow{2}{4em}{Prop.~\ref{prop:harmonic}}\\
  & $(\ln K)M_{r,K}$ & (asymptotic formula) &&  \\\hline
  $r\in(-\infty,-1)$ & $\frac{r}{r+1}K^{1+1/r}M_{r,K}$ & asymptotically precise & &  Prop.~\ref{prop:non-integrable-optimality}\\\hline
  $r=-\infty$ & $K M_{r,K}$ & precise & Bonferroni&\\
  \hline
\end{tabular}
\end{table}

Two of our results, Theorems~\ref{thm:Bernard} and~\ref{thm:Embrechts}, are general
in the sense of describing properties of Kolmogorov's general averaging function~\eqref{eq:mean};
Theorem~\ref{thm:Bernard} covers the case of integrable $\phi$
(and therefore, the case $r>-1$ in \eqref{eq:mean-r}),
and Theorem~\ref{thm:Embrechts} covers the case of non-integrable $\phi$
(and therefore, the case $r\le-1$ in \eqref{eq:mean-r}).
Propositions \ref{prop:asymptotically-precise}-\ref{prop:harmonic} 
give results shown in Table \ref{tab:main}.
We further show in  Proposition \ref{prop:homo} that, for the definition of a p-variable in \eqref{eq:p},
one can replace ``for all $\epsilon$'' by ``for some $\epsilon$'' in all averaging methods
based on \eqref{eq:mean-r}.

\begin{theorem}\label{thm:Bernard}
  Suppose a continuous strictly monotonic $\phi:[0,1]\to[-\infty,\infty]$ is integrable,
  i.e., $\int_0^1\left|\phi(u)\right|\dd u<\infty$.
  Then, for any $K\in\{2,3,\dots\}$ and any $\epsilon>0$,
  \begin{equation}\label{eq:general}
    \Prob
    \left(
      M_{\phi,K}(p_1,\dots,p_K)
      \le
      \psi
      \left(
        \frac{1}{\epsilon}
        \int_0^{\epsilon}
          \phi(u)
        \dd u
      \right)
    \right)
    \le
    \epsilon.
  \end{equation}
\end{theorem}

As we stated it, Theorem~\ref{thm:Bernard} gives a critical region of size $\epsilon$.
An alternative statement is that $\Psi^{-1}(M_{\phi,K})$ is a merging function,
where the strictly increasing function $\Psi$ is defined by
\begin{equation}\label{eq:term}
  \Psi(\epsilon)
  :=
  \psi
  \left(
    \frac{1}{\epsilon}
    \int_0^{\epsilon}
      \phi(u)
    \dd u
  \right),
  \quad
  \epsilon\in(0,1).
\end{equation}

In what follows, the expression $(r+1)^{1/r}$ is understood to be $e=\lim_{r\to0}(r+1)^{1/r}$ when $r=0$
and $1=\lim_{r\to\infty}(r+1)^{1/r}$ when $r=\infty$.

Using Theorem~\ref{thm:Bernard} with $\phi(u)=u^r$,
we see that   $
    (r+1)^{1/r}
    M_{r,K}
  $
  is a merging function for $r>-1$. 
Moreover, we can show the constant $(r+1)^{1/r}$  cannot be improved in general.

\begin{proposition}\label{prop:asymptotically-precise}
  Let $r\in (-1,\infty]$.
  Then 
  $
    (r+1)^{1/r}
    M_{r,K}
  $, $K=2,3,\dots$,
  is   a  family   of merging functions and it is asymptotically precise.
\end{proposition}

In particular, Proposition \ref{prop:asymptotically-precise} implies
that for the geometric mean, it has to be multiplied by the constant factor $e$,
which cannot be improved in general for large $K$.
The next proposition characterizes the condition
under which the merging function in Proposition \ref{prop:asymptotically-precise} is precise.

\begin{proposition}\label{prop:precise}
  For  $r\in(-1,\infty)$ and $K\in\{2,3,\dots\}$,
  the merging function $ (r+1)^{1/r}M_{r,K}$
  is precise if and only if $r\in[\frac{1}{K-1},K-1]$.
\end{proposition}

The most straightforward yet relevant example of Proposition~\ref{prop:precise},
the arithmetic average multiplied by $2$,
namely,
\[
  M_{1,K}(p_1,\dots,p_K)
  :=
  \frac{2}{K}\sum_{k=1}^K p_k,
\]
is a precise merging function for all $K\ge 2$, a result obtained in \cite{Ruschendorf:1982}. 
As another special case of Proposition~\ref{prop:precise},
the scaled quadratic average multiplied by $\sqrt 3$, namely $\sqrt{3} M_{2,K}$, is a merging function,
and it is precise if and only if $K\ge3$.

In the case $r\ge1$, the merging function in Proposition~\ref{prop:precise}
can be modified in an explicit way such that it remains precise even for $r>K-1$.

\begin{proposition}\label{prop:super-precise}
  For $K\in\{2,3,\dots\}$ and $r\in[1,\infty)$, the function
  $
    \min(r+1,K)^{1/r}
    M_{r,K}
 $
  is a precise merging function.
\end{proposition}

Because of the importance of geometric mean as one of the Platonic means,
the following result gives a precise (albeit somewhat implicit) expression
for the corresponding precise merging function.

\begin{proposition}\label{prop:geometric}
  For each $K\in\{2,3,\dots\}$,
  $a^G_K M_{0,K}$ is a precise merging function,
  where
  \[
    a^G_K
    :=
    \frac{1}{c_K}
    \exp
    \left(
      -(K-1)(1-K c_K)
    \right)
  \]
  and $c_K$ is the unique solution to the equation
  \begin{equation}\label{eq:ga}
    \ln(1/c-(K-1))
    =
    K - K^2 c
  \end{equation}
  over $c\in(0,1/K)$. Moreover, $a^G_K\le e$ and  $a^G_K\to e$ as $K\to \infty$.
\end{proposition}

Table~\ref{tab:geometric} reports several values of $a^G_K/e$ calculated numerically
and suggests that in practice there is no point in improving the factor $e$ for $K\ge5$.

\begin{table}[htbp]
\caption{Numeric values of $a^G_K/e$ for the geometric mean}
\label{tab:geometric}
\begin{center}
\begin{tabular}{cc|cc|cc}
  \hline
  $K$ & $a^G_K/e$ & $K$ & $a^G_K/e$ & $K$ & $a^G_K/e$\\
  \hline
  2 & 0.7357589 & 5 & 0.9925858 & 10 & 0.9999545\\
  3 & 0.9286392 & 6 & 0.9974005 & 15 & 0.9999997\\
  4 & 0.9779033 & 7 & 0.9990669 & 20 & 1.0000000\\
  \hline
\end{tabular}
\end{center}
\end{table}

The condition $r>-1$ in Proposition \ref{prop:asymptotically-precise}  ensures that the term~\eqref{eq:term} is finite,
and also that the condition $\int_0^1\left|\phi(u)\right|\dd u<\infty$ in Theorem~\ref{thm:Bernard} is satisfied.
However, the condition rules out the harmonic mean (for which $r=-1$) and the minimum ($r=-\infty$).
The next simple corollary of another known result
(Theorem~4.2 of \cite{Embrechts/Puccetti:2006};
see also Theorem 2.3 of \cite{Puccetti2013sharp})
will cover these cases as well.

\begin{theorem}\label{thm:Embrechts}
  Suppose $\phi:[0,1]\to[-\infty,\infty]$ is a strictly decreasing continuous function satisfying $\phi(0)=\infty$.
  Then, for any $\epsilon\in(0,1)$ such that $\phi(\epsilon)\ge0$,
  \begin{equation}\label{eq:Embrechts}
    \Prob
    \left(
      M_{\phi,K}(p_1,\dots,p_K)
      \le
      \epsilon
    \right)
    \le
    \inf_{t\in(0,\phi(\epsilon)]}
    \frac{\int_{\phi(\epsilon)-t}^{\phi(\epsilon)+(K-1)t}\psi(u)\dd u}{t}.
  \end{equation}
\end{theorem}

As $t\to0$, the upper bound in~\eqref{eq:Embrechts} is not informative since, for $t\approx0$,
\[
  \frac{\int_{\phi(\epsilon)-t}^{\phi(\epsilon)+(K-1)t}\psi(u)\dd u}{t}
  \approx
  \frac{Kt\psi(\phi(\epsilon))}{t}
  =
  K\epsilon,
\]
which is dominated by the Bonferroni bound.
On the other hand, the upper bound is informative when $t=\phi(\epsilon)$
provided the integral is convergent.
For example, we can see that for $r<-1$,
$
  \frac{r}{r+1}
  K^{1+1/r}
  M_{r,K}
$
is a merging function.
In what follows, the term $r/(r+1)$ should be understood as its limit $1$ when $r=-\infty$.

\begin{proposition}\label{prop:non-integrable-optimality}
  Let $r\in[-\infty,-1)$.
 Then 
  $
    \frac{r}{r+1}
    K^{1+1/r}
    M_{r,K}
  $, $K=2,3,\dots$,
  is   a  family   of merging functions and it is asymptotically precise. 
\end{proposition}

Proposition \ref{prop:non-integrable-optimality} includes the Bonferroni bound~\eqref{eq:Bonferroni} as special case:
for $r:=-\infty$, we obtain that $K M_{-\infty,K}$ is a merging function. 
On the other hand, Proposition \ref{prop:non-integrable-optimality}  does not cover the case $r=-1$ of harmonic mean directly,
but easily implies a bound  leading to the merging function $e \ln K M_{-1,K}$, which turns out to be not so crude. 

\begin{proposition}\label{prop:harmonic}
  Set $a^{H}_K := \frac{(y_K+K)^2}{(y_K+1)K}$, $K>2$,
  where $y_K$ is the unique solution to the equation
  \begin{equation*}
    y^2 = K((y+1)\ln (y+1)-y),
    \qquad
    y\in(0,\infty).
  \end{equation*}
  Then $a^{H}_K  M_{-1,K}$ is a precise merging function.
  Moreover, $a^{H}_K  \le e\ln K$ and $a^{H}_K/\ln K \to 1$ as $K\to\infty$.
\end{proposition}

Even though $a^{H}_K/\ln K\to 1$, the rate of convergence is very slow,
and $a^{H}_K/\ln K > 1$ for moderate values of $K$.
In practice,
it might be better to use the conservative merging function $(e\ln K) M_{-1,K}$. 
Table \ref{tab:harmonic} reports several values of $a^{H}_K/\ln K$ calculated numerically.
For instance, for $K \ge 10$, one may use $(2\ln K) M_{-1,K}$,
and for $K\ge 50$, one may use $(1.7 \ln K)M_{-1,K}$.

\begin{table}[htbp]
\caption{Numeric values of $a^H_K/\ln K$ for the harmonic mean}
\label{tab:harmonic}
\begin{center}
\begin{tabular}{cc|cc|cc}
  \hline
  $K$ & $a^H_K/\ln K$ & $K$ & $a^H_K/\ln K$ & $K$ & $a^H_K/\ln K$\\
  \hline
  3 & 2.499192 & 10 & 1.980287 & 100 & 1.619631\\
  4 & 2.321831 & 20 & 1.828861 & 200 & 1.561359\\
  5 & 2.214749 & 50 & 1.693497 & 400 & 1.514096\\
  \hline
\end{tabular}
\end{center}
\end{table}

The main emphasis of this section has been on characterizing $a>0$ such that $F:=a M_{r,K}$ is a merging function,
or a precise merging function.
Recall that $F:[0,1]^K\to[0,\infty)$ is a merging function
if and only if \eqref{eq:merging} holds for all $\epsilon\in(0,1)$,
and that $F$ is a precise merging function
if and only if \eqref{eq:precise-merging} holds for all $\epsilon\in(0,1)$.
The next proposition shows that in both statements ``for all'' can be replaced by ``for some''
if $F=a M_{r,K}$.
A practical implication is that even if an applied statistician is interested in the property of validity \eqref{eq:p}
only for specific values of $\epsilon$ (such as $0.01$ or $0.05$)
and would like to use $a M_{r,K}$ as a merging function,
she is still forced to ensure that \eqref{eq:p} folds for all $\epsilon$.

\begin{proposition}\label{prop:homo}
  For any $a>0$, $r\in[-\infty,\infty]$, and $K\in\{2,3,\dots\}$:
  \begin{itemize}
  \item[(a)]
    $F:=a M_{r,K}$ is a merging function if and only if~\eqref{eq:merging}
    holds for some $\epsilon\in(0,1)$;
  \item[(b)]
    $F:=a M_{r,K}$ is a precise merging function
    if and only if \eqref{eq:precise-merging} holds for some $\epsilon\in(0,1)$.
  \end{itemize}
\end{proposition}

\section{Combining p-values by weighted averaging}
\label{sec:weighted}

In this section we will briefly consider a more general notion of averaging:
\begin{equation*}
  M_{\phi,\w}(p_1,\dots,p_K)
  :=
  \psi
  \left(
    w_1\phi(p_1)+\dots+w_K\phi(p_K)
  \right)
\end{equation*}
in the notation of~\eqref{eq:mean},
where $\w=(w_1,\dots,w_K)\in\Delta_K$ is an element of the standard $K$-simplex
\[
  \Delta_K
  :=
  \left\{
    (w_1,\dots,w_K)\in[0,1]^K
    \st
    w_1+\dots+w_K=1
  \right\}.
\]
One might want to use a weighted average in a situation
where some of p-values are based, e.g., on bigger experiments,
and then we might want to take them with bigger weights.
Intuitively, the weights reflect the prior importance of the p-values
(see, e.g., \cite[p.~5]{Borenstein/etal:2009} for further details).

Much fewer mathematical results in the literature are available for asymmetric risk aggregation. 
For this reason, we will concentrate on the easier integrable case, namely, $r>-1$.
Theorem~\ref{thm:Bernard} can be generalized as follows.

\begin{CustomTheorem}{\ref{thm:Bernard}w}\label{thm:Bernard-w}
  Suppose a continuous strictly monotonic $\phi:[0,1]\to[-\infty,\infty]$ is integrable and $\w\in\Delta_K$.
  Then, for any $\epsilon>0$,
  \begin{equation*}
    \Prob
    \left(
      M_{\phi,\w}(p_1,\dots,p_K)
      \le
      \psi
      \left(
        \frac{1}{\epsilon}
        \int_0^{\epsilon}
          \phi(u)
        \dd u
      \right)
    \right)
    \le
    \epsilon.
  \end{equation*}
\end{CustomTheorem}

Similarly to \eqref{eq:mean-r}, we set
\begin{equation*}
  M_{r,\w}(p_1,\dots,p_K)
  :=
  \left(
    w_1 p_1^r + \dots + w_K p_K^r
  \right)^{1/r}
\end{equation*}
for $r\in\R$ and $\w=(w_1,\dots,w_K)\in\Delta_K$.
We can see that Proposition \ref{prop:asymptotically-precise} still holds when $M_{r,K}$ is replaced by $M_{r,\w}$,
for any $r\in\R$ and $\w\in\Delta_K$.
This is complemented by the following proposition, which is the weighted version of Proposition \ref{prop:precise}.

\begin{CustomProposition}{\ref{prop:precise}w}\label{prop:precise-w}
  For $\w=(w_1,\dots,w_K)\in\Delta_K$ and $r\in(-1,\infty)$,
  the merging function
  $(r+1)^{1/r}M_{r,\w}$
  is precise if and only if $w\le1/2$ and $r\in[\frac{w}{1-w},\frac{1-w}{w}]$,
  where $w:=\max_{k=1,\dots,K}w_k$.
\end{CustomProposition}

Next we generalize Proposition~\ref{prop:super-precise} to non-uniform weights.

\begin{CustomProposition}{\ref{prop:super-precise}w}\label{prop:super-precise-w}
  For $\w=(w_1,\dots,w_K)\in\Delta_K$ and $r\in[1,\infty)$,
  the function
  $
    \min(r+1,\frac{1}{w})^{1/r}
    M_{r,\w}
  $
  is a precise merging function,
  where $w:=\max_{k=1,\dots,K}w_k$.
\end{CustomProposition}

An interesting special case of Proposition~\ref{prop:super-precise-w} is for $r=1$ (weighted arithmetic mean).
If $w\le1/2$, i.e., no single experiment outweighs the total of all the other experiments,
the optimal multiplier for the weighted average is 2, exactly as in the case of the arithmetic average.
If $w>1/2$, i.e., there is a single experiment that outweighs all the other experiments,
our merging function is, assuming $w_1=w$,
\[
  \frac 1 w M_{1,\w}(p_1,\dots,p_K)
  =
  p_1
  +
  \sum_{k=2}^K
  \frac{w_k}{w}
  p_k.
\]
It is obtained by adding weighted adjustments to the p-value obtained from the most important experiment.

\section{Efficiency of merging functions}
\label{sec:efficiency}

So far we have emphasized the \emph{validity} of our methods of combining p-values:
the combined p-value $P$ is guaranteed to satisfy \eqref{eq:p} under the null hypothesis
$p_1,p_2,\dots,p_K\in\mathcal{U}$.
In this section we will discuss their \emph{efficiency}:
under alternative hypotheses, we would like the combined p-value to be small.

We will be interested in asymptotic results as $K\to\infty$.
Recall that our merging functions are
\[
  P_{r,K}(p_1,\dots,p_K)
  :=
  a_{r,K}
  M_{r,K}(p_1,\dots,p_K),
\]
where $r\in[-\infty,\infty]$, $P_{r,K}$ is the combined p-value, and
\begin{equation*}
  a_{r,K}
  :=
  \begin{cases}
    (r+1)^{1/r} & \text{if $r\in(-1,\infty]$}\\
    e\ln K & \text{if $r=-1$}\\
    \frac{r}{r+1}K^{1+1/r} & \text{if $r\in[-\infty,-1)$},
  \end{cases}
\end{equation*}
where for $r=0$, $r=\infty$ and $r=-\infty$ one uses the limiting values of $e$, $1$, and $K$, respectively.
Notice that we do not truncate $P_{r,K}$ by replacing it with $\min(P_{r,K},1)$
and interpret large values of $P_{r,K}$, $P_{r,K}\gg1$,
as indicators of the weakness of the merging function:
they show us how far it is from being useful.

The rest of this section consists of two parts.
In Subsection~\ref{subsec:IID} we consider an infinite sequence of IID p-values $p_1,p_2,\dots$
and analyze the performance of our merging functions on $p_1,\dots,p_K$ as $K\to\infty$.
Of course, we are not interested in independent $p_1,p_2,\dots$ per se;
after all, if we know the p-values to be independent,
we should use much more efficient methods, such as Fisher's, that assume independence.
The result in Subsection~\ref{subsec:IID} will serve as a basic tool
for the analysis of the  general symmetric framework, considered in Subsection~\ref{subsec:exchangeability}.
We assume only the exchangeability of $p_1,p_2,\dots$,
and according to de Finetti's theorem every exchangeable probability measure on $[0,1]^{\infty}$
is a mixture of IID components
(in particular, when we look at the realized sequence of p-values,
we can only draw conclusions about the realized IID component).

\subsection{IID p-values}
\label{subsec:IID}

\begin{figure}[bt]
  \begin{center}
    \includegraphics[width=0.5\textwidth]{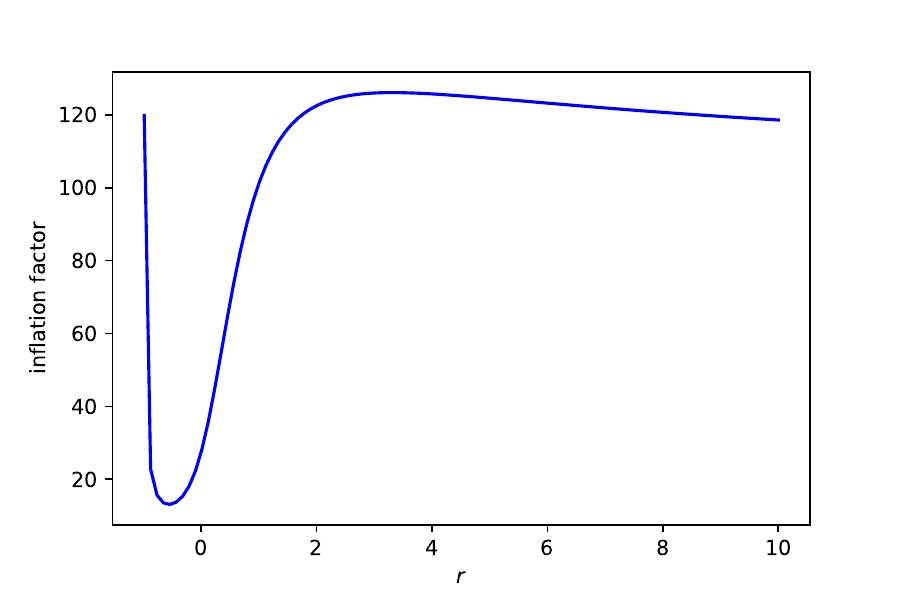}
  \end{center}
  \caption{The asymptotic inflation factor $\lim_{K\to\infty}P_{r,K}/p$
    as function of $r$ when each p-value $p_k$ is either $p$ or $100p$ (with equal probabilities),
    for some $p\in(0,0.01)$.
    (The trivial factor is 100, and the ideal, but unattainable, factor is 1.)}
  \label{fig:c05}
\end{figure}

In this subsection we assume that the p-values $p_1,p_2,\dots$
are generated independently from the same probability measure $Q$ on $[0,1]$.
If the support of $Q$ is finite and does not include 0,
the asymptotic performance of $P_{r,K}$ is poor for $r<-1$,
but some of $r>-1$ may have reasonable performance:
see Figure~\ref{fig:c05} for an example.

In general we impose integrability conditions on $p_k$.
Set
\[
  \Pi
  =
  \Pi(Q)
  :=
  \sup
  \left\{
    m\in[0,\infty) \mid \Expect(p_1^{-m})=\int p^{-m}Q(\dd p)<\infty
  \right\}
  \in
  [0,\infty].
\]
By the Lyapunov inequality \cite[Theorem~3.2.5]{Gut:2013},
$\Expect(p_1^{-m})<\infty$ for $m<\Pi$
and $\Expect(p_1^{-m})=\infty$ for $m>\Pi$.
The value of $\Pi$ reflects the asymptotic power of the p-values;
we will discuss this further at the end of this subsection.

\begin{figure}[bt]
  \begin{center}
    \includegraphics[width=0.6\textwidth,trim={0mm 35mm 0mm 35mm},clip]{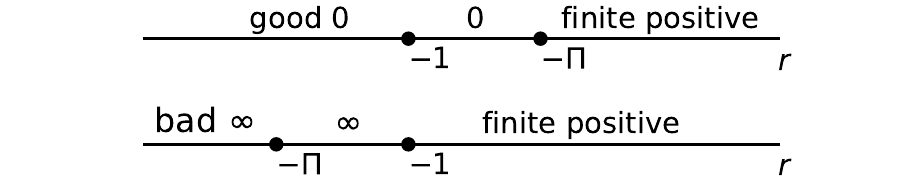}
  \end{center}
  \caption{The combined p-value for different $r$
  in the cases $\Pi<1$ (top) and $\Pi>1$ (bottom).}
  \label{fig:summary}
\end{figure}

The following proposition describes the performance of our merging functions
for various values of $r$ (except for the two critical values, $-1$ and $-\Pi$).
To analyze the asymptotic behaviour of $P_{r,K}$ as $K\to\infty$,
we will write $P_{r,K}\approx K^a$ for $a\in\R$ if
\begin{equation}\label{eq:approx}
   \liminf_{K\to\infty}
  \frac{P_{r,K}}{K^b}
  =     0   \text{~~for $b>a$ ~and ~} 
  \lim_{K\to\infty}
  \frac{P_{r,K}}{K^b}
  =    \infty  \text{~~for $b<a$}
  \quad
  \text{a.s.} 
\end{equation} 
This means that $P_{r,K}$ is roughly of the order $K^a$, as least in the sense of $\liminf$.
Also recall that for two sequences $b_K$ and $c_K$, $K=1,2,\dots$, the symbol $b_K\sim c_K$ means $b_K/c_K\to 1$.

\begin{proposition}\label{prop:IID}
  For $\Pi>1$, we have:
  \begin{enumerate}
  \item[(i)]
    If $r>-1$, then
    \begin{equation}\label{eq:big}
      P_{r,K}
      \to
      (r+1)^{1/r}
      \Expect(p_1^r)^{1/r}
      \in
      (0,\infty)
      \quad
      \text{a.s.}
    \end{equation}
  \item[(ii)]
    If $-\Pi<r<-1$, then
    \begin{equation*}
      P_{r,K}
      \sim
      \frac{r}{r+1}
      K^{1+1/r}
      \Expect(p_1^r)^{1/r}
      \to
      \infty
      \quad
      \text{a.s.}
    \end{equation*}
  \item[(iii)]
    If $-\infty<r<-\Pi$, then $P_{r,K}\approx K^{1-1/\Pi} \to \infty$.
  \item[(iv)]
    If $r=-\infty$, then $\lim _{K\to\infty}  {P_{r,K}/ K^{b}} = \infty$ in probability for $b< {1-1/\Pi} $. 
  \end{enumerate}
  For $\Pi\in (0,1)$, we have:
  \begin{enumerate}
  \item[(i)]
    If $r>-\Pi$, then \eqref{eq:big} holds.
  \item[(ii)]
    If $-1<r<-\Pi $, then $P_{r,K}\approx K^{-1/r-1/\Pi}\to 0$.
  \item[(iii)]
    If $-\infty<r<-1$, then $P_{r,K}\approx K^{1-1/\Pi}\to 0$.
  \item[(iv)]
    If $r=-\infty$, then $\liminf_{K\to\infty} P_{r,K}/K^b = 0$ a.s.  for $b> {1-1/\Pi} $.  
  \end{enumerate}
\end{proposition}

The results given in Proposition~\ref{prop:IID} are roughly summarized
in Figure~\ref{fig:summary}, where $r\in\{-\infty,\infty\}$ is also allowed.
``Finite positive'' means that $P_{r,K}$ converges to a finite positive number as $K\to\infty$,
$0$ means that it converges to $0$ (at least in the sense of $\liminf$),
``good $0$'' means that the rate of convergence to $0$ is optimal
(within the accuracy of \eqref{eq:approx}),
$\infty$ means that it converges to $\infty$,
and ``bad $\infty$'' means the fastest rate of convergence to $\infty$. 

As we mentioned earlier, in the IID situation, the natural interpretation of $\Pi$
is that it measures the asymptotic power of the p-values.
When $p_1$ is uniformly distributed on $[0,1]$ (as it is under the null hypothesis),
$\Pi=1$.
If $\Pi>1$, the p-values are asymptotically powerless for rejecting the null hypotheses,
and when $\Pi<1$, there is some power
(e.g., the density function of $p_1$ is $f(p):=cp^{c-1}$ for $c\in(0,1)$,
as in, e.g., \cite{Sellke/etal:2001}).
The general message of this subsection is that, in the IID case, using the Bonferroni merging function is a safe option,
at least asymptotically,
since it performs well in the interesting case $\Pi<1$ (top of Figure~\ref{fig:summary});
it might be also interesting to note that any merging function with $r<-1$ would achieve the same best rate of convergence to $0$.
The situation becomes more complicated in the case of dependent p-values,
which will be considered in the next subsection.

\subsection{Exchangeable p-values}
\label{subsec:exchangeability}

Now we relax our IID assumption to the more interesting case that the sequence $p_1,p_2,\dots$ is exchangeable.
This most general symmetry assumption includes cases of heavy dependence.
By de Finetti's theorem,
every exchangeable distribution on $[0,1]^{\infty}$ is a mixture
$\int \tilde Q ^{\infty}\nu(\dd\tilde Q)$ of IID distributions $\tilde Q^\infty$;
see, e.g., \cite[Theorem 1.49]{Schervish:1995}.
We can interpret $\nu$ as the limiting empirical distribution function of $p_1,\dots,p_K$.

The performance of our merging functions depends on the distribution of $\Pi(\tilde Q)$.
Applying Proposition~\ref{prop:IID} in combination with de Finetti's theorem,
we obtain the following two cases:
\begin{enumerate}
\item
  If $\nu(\Pi(\tilde Q)<1)=1$,
  the values $r\in[-\infty,-1)$ (including Bonferroni)
  perform very well
  (they are in the ``good $0$'' area in Figure \ref{fig:summary}).
\item
  If $\nu(\Pi(\tilde Q)>1)>0$,
  the values $r\in[-\infty,-1)$ lead to an infinite expected combined p-value 
  (or, with a positive probability, to a combined p-value of 1
  when truncated) as $K\to\infty$.
\end{enumerate}

Recall that the main motivation for our merging methods is the existence of unknown (possibly heavy) dependence.
With dependence among p-values, the case $\nu(\Pi(\tilde Q)>1)>0$ is not unusual.
For instance, in the extreme case of perfect dependence $p_1=\dots=p_K$,
if $\Prob(p_1=0)=0$, then $\nu(\Pi(\tilde Q)=\infty)=1$.

To analyze the influence of the strength of dependence on the performance of $r\in [-\infty,\infty]$,
we consider the following simple example.  
Let
\begin{equation}\label{eq:defx} 
  p_k = \Phi(X_k),
  \quad
  X_k= \rho Z + \sqrt{1-\rho^2} Z_k - \mu,
  \quad
  k=1,\dots,K,
\end{equation}
where $\Phi$ is the standard normal distribution function,
$Z,Z_1,\dots,Z_K$ are IID standard  normal random variables,
and $\mu\ge 0$ and $\rho\in [0,1]$ are constants.
In other words, $p_k$ is the p-value resulting from the $k$th standard one-sided z-test
of the null hypothesis $\mu=0$ against the alternative $\mu>0$
using the statistic $X_k$ from $\mathrm{N}(-\mu,1)$ with unknown $\mu$, $k=1,\dots,K$.
Note that $\rho=0$ corresponds to the case where $p_1,\dots,p_K$ are independent,
while $\rho=1$ corresponds to the case where $p_1,\dots,p_K$ are perfectly dependent.

Obviously, the model of $(p_1,\dots,p_K)$ is exchangeable,
and the marginal distribution $Q$ of $(p_1,\dots,p_K)$ does not depend on the correlation $\rho$.
For $z\in\R$, let $Q_z$ be the distribution of $\Phi(\sqrt{1-\rho^2 }W + \rho z - \mu)$,
where $W$ is a standard normal random variable.
Clearly, conditional on $Z=z$, $p_1,\dots,p_K$ are IID with distribution $Q_z$. 
Note that $\Pi(Q_z)=1/(1-\rho^2)>1$ if $\rho>0$
($\Pi(Q_z)$ is computed in Lemma \ref{lem:gaussian} on p.~\pageref{lem:gaussian} in Appendix~\ref{app:proofs-validity}).
Therefore, in the presence of positive dependence,
we are in Case~2 above ($\nu(\Pi( Q_z)>1)=1$); thus some choice of $r\ge -1$ may be optimal asymptotically as $K\to \infty$. 
For finite $K$, some simulation results of $P_{r,K}$ for the model \eqref{eq:defx} are reported in Subsection \ref{sec:74}.

\section{Choosing the merging function}
\label{sec:choice}

In this section we discuss, for simplicity  only in the symmetric case, how to  choose  a merging function. 

\subsection{A rule of thumb}

First we state a crude rule of thumb for choosing $r$.
Since any method based on the observed values of $p_1,\dots,p_K$ would affect the validity of the method
(see Subsection \ref{subsec:search}),
we have to rely on prior or side information for a suitable choice of $r$.
As a rule of thumb, if there is potentially substantial dependence among the p-values,
then we should not use Bonferroni, and the harmonic mean might be a safer choice.
If we are certain that the dependence is really strong,
then the geometric and the arithmetic means might be an even better option.
See Subsection \ref{sec:74} for a simulation study illustrating this point.

\subsection{Practical issues}
\label{subsec:practical}

A fairly wide family of merging functions is provided in this paper.
For a practitioner, given access to a variety of merging functions,
it might be tempting to try many of them
and then pick one of the merging functions that work well
(perhaps the one yielding the smallest p-value).
Of course, this is not a valid approach, and its inadmissibility has been discussed at length in recent literature
(see, e.g., \cite{AS:2019}).
A search over various merging functions itself needs a multiple-testing correction,
and the simplest procedure in the context of this paper
would be to apply another merging function in our family
to the outputs of the merging functions that we used in our search.
This will be the topic of Subsection~\ref{subsec:search}.

It might be also tempting to divide the available p-values into two parts,
to find a suitable value of $r$ from one part,
and then compute the combined p-value using $a_{r,K}M_{r,K}$ applied to the other part.
However, this is not a valid approach if dependence is present among the p-values.
Under possibly heavy dependence,
any peeking into the data is likely to destroy validity
(e.g., learning any of the p-values makes the conditional distribution of the remaining p-values non-uniform).
From the vantage point of Section~\ref{sec:efficiency},
the fact that we can never observe any IID components apart from the realized one
prevents us from using interactive methods of adaptation such as the STAR method in \cite{Lei/etal:2017}.

\subsection{Combining merging functions}
\label{subsec:search}

In this subsection we will be interested in valid,
necessarily very restrictive,
ways of searching for a good merging function.
Namely, we consider the following mechanism:
for fixed $r_1,\dots,r_m\in [-\infty,\infty]$ and partition $(A_1,\dots,A_m)$ of $[0,1]^K$,
we would like to use
\begin{equation}\label{eq:adaptive}
  F(p_1,\dots,p_K)
  :=
  b
  \sum_{i=1}^m
  P_{r_i,K} (p_1,\dots,p_K)
  \III_{A_i}(p_1,\dots,p_K)
\end{equation}
as a merging function,
where $b$ is a positive constant. 
In this case the observed p-values $p_1,\dots,p_K$
determine which of the merging functions in our family $P_{r,K}$ is used,
and the partition $(A_1,\dots,A_m)$ serves as decision criterion.
The value of $b$ is chosen in such a way that $F$ is a valid merging function.

A natural and simple choice of the decision criterion in \eqref{eq:adaptive}
is to take the minimum of the merging functions $P_{r,K}$ to exploit their power.
That is, by choosing $A_i$ as the region in $[0,1]^K$ where $P_{r_i,K}(p_1,\dots,p_K)$ is the smallest,
we arrive at
\begin{equation}\label{eq:adaptive2}
  F(p_1,\dots,p_K)
  :=
  b
  \min_{i=1,\dots,m}
  P_{r_i,K} (p_1,\dots,p_K).
\end{equation}
The constant $b$ in \eqref{eq:adaptive2} may be seen as the price to pay
to exploit the power of different merging functions,
and it is typically larger than $1$.
Clearly, choosing $b=m$ leads to a valid merging function in \eqref{eq:adaptive2},
and this is precisely applying the Bonferroni method on the combined p-values $P_{r_1,K},\dots,P_{r_m,K}$.
Simple examples show that in some cases a smaller value of $b$ will also lead to a valid merging function:
e.g., if all $r_1,\dots,r_m$ are similar,
then $b$ can be chosen close to $1$.

We consider in detail the most applicable cases where we search over only two values of $r$,
one of them being $-\infty$ (the Bonferroni merging function).
We would like to be competitive with the best of the two values,
and so combine the two merging functions using the Bonferroni merging function.
Namely, we consider the compound \emph{Bonferroni-arithmetic} (BA) merging function
\begin{equation}\label{eq:BA}
  F_K^{\BA}(p_1,\dots,p_K)
  :=
  2
  \min
  \left(
    K\min(p_1,\dots,p_K),
    2\bar p
  \right) 
\end{equation}
where $\bar p$ is the arithmetic mean of $p_1,\dots,p_K$,
and  the compound \emph{Bonferroni-geometric} (BG) merging function
\begin{equation}\label{eq:BG}
  F_K^{\BG}(p_1,\dots,p_K)
  :=
  2
  \min
  \left(
    K\min(p_1,\dots,p_K),
    e\tilde p
  \right) 
\end{equation} 
where $\tilde p$ is the geometric mean of $p_1,\dots,p_K$.
Obviously, both methods are valid merging functions.
Moreover, it turns out that these two merging functions are asymptotically precise.
In other words,
the price to pay for exploiting the power of both the Bonferroni method and the arithmetic/geometric average
is precisely a factor of $2$.
\begin{proposition}\label{prop:BA-BG}
  Both families of merging functions $F_K^{\BA}$ and $F_K^{\BG}$, $K=2,3,\dots$, in \eqref{eq:BA} and \eqref{eq:BG} are asymptotically precise.   
\end{proposition}
A crucial advantage of $F_K^{\BA}$ and $F_K^{\BG}$ over the merging functions in our family $a_{r,K}M_{r,K}$
is that both $F_K^{\BA}$ and $F_K^{\BG}$ improve greatly on the performance
crudely shown in Figure \ref{fig:summary}:
for $F_K^{\BA}$ and $F_K^{\BG}$,
in the top plot, we will get ``finite positive'' (the best value in that plot),
and in the bottom plot, we will get ``good 0'' (also the best value in that plot).
In particular, neither $F_K^{\BA}$ nor $F_K^{\BG}$ is dominated by $a_{r,K}M_{r,K}$ for any $r$.

\subsection{A simulation study}\label{sec:74}

\begin{figure}[htbp]
  \begin{center}
    \includegraphics[width=0.48\textwidth]{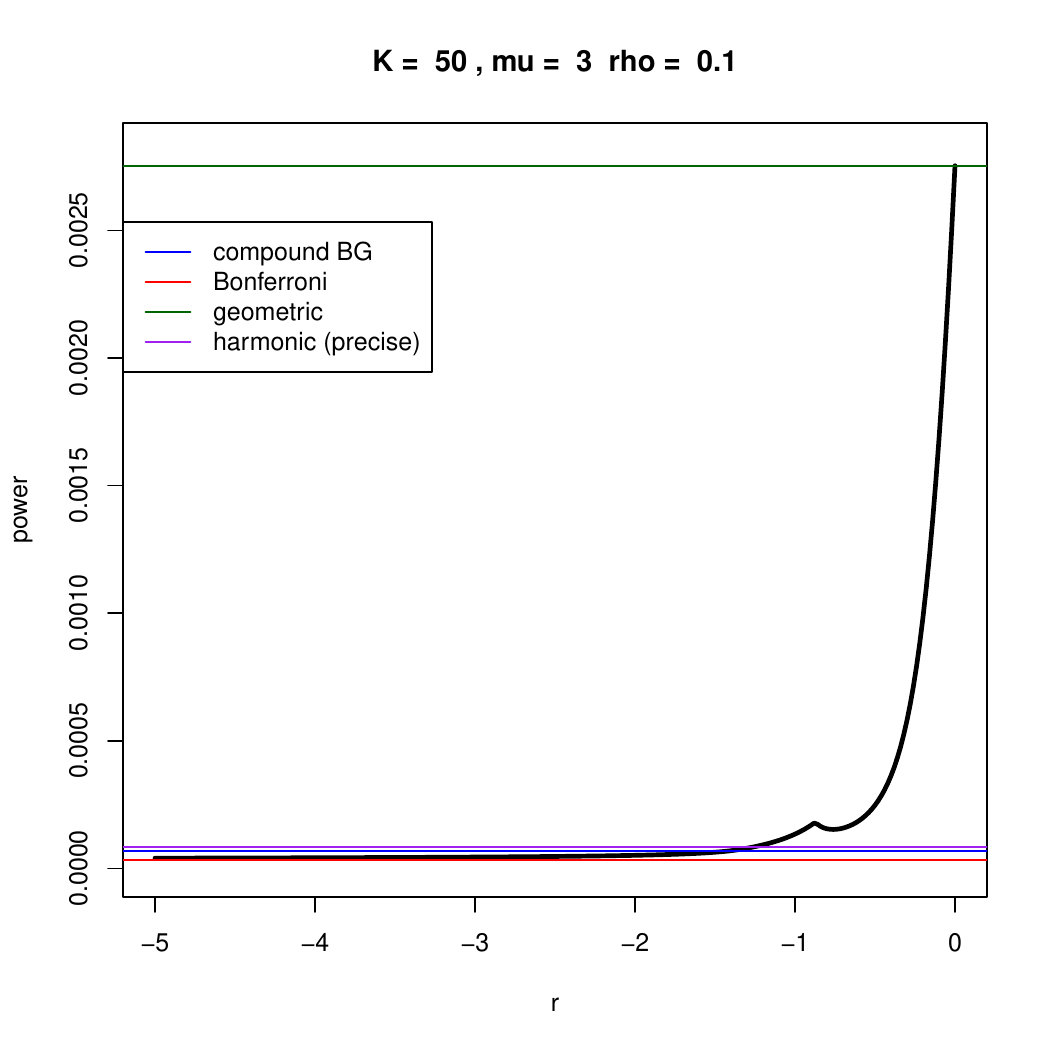}
    \includegraphics[width=0.48\textwidth]{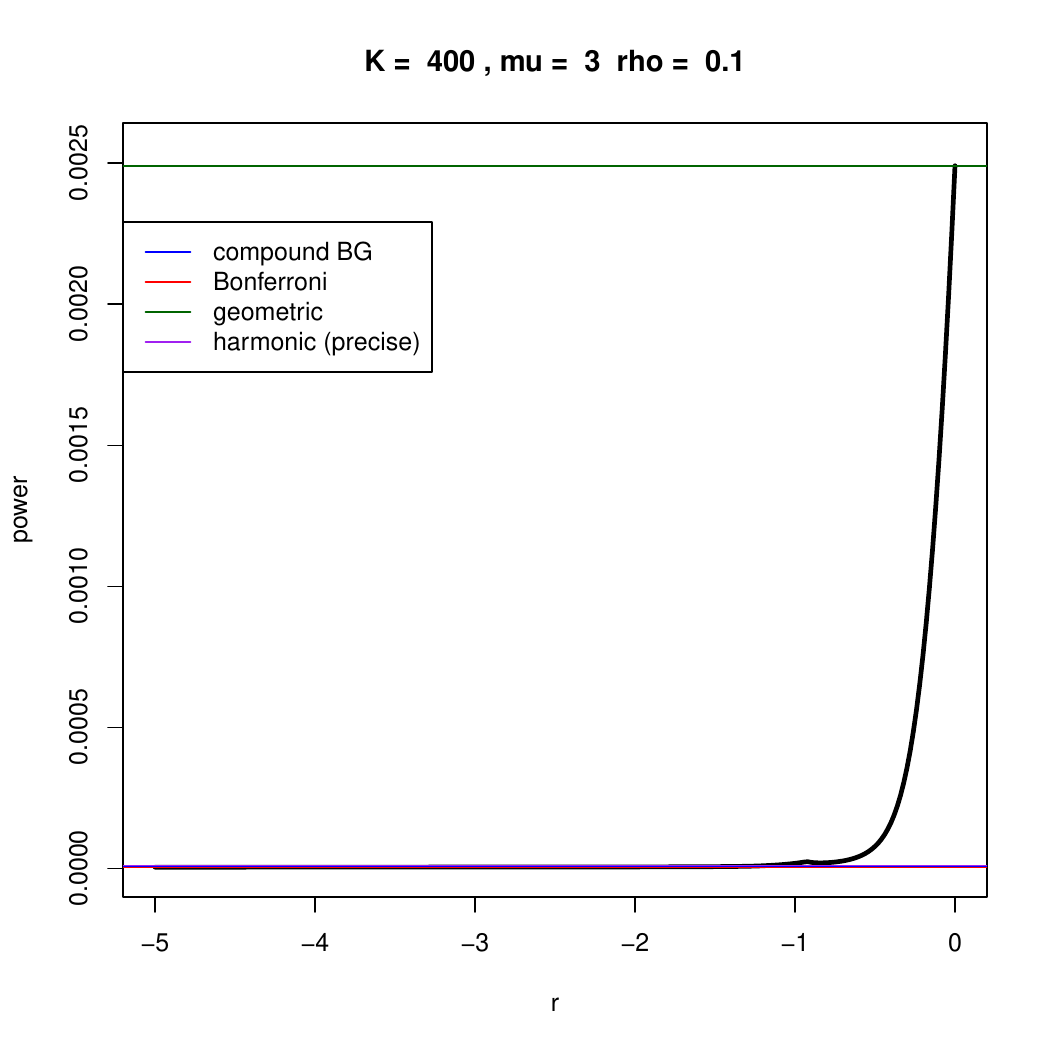}\\
    \includegraphics[width=0.48\textwidth]{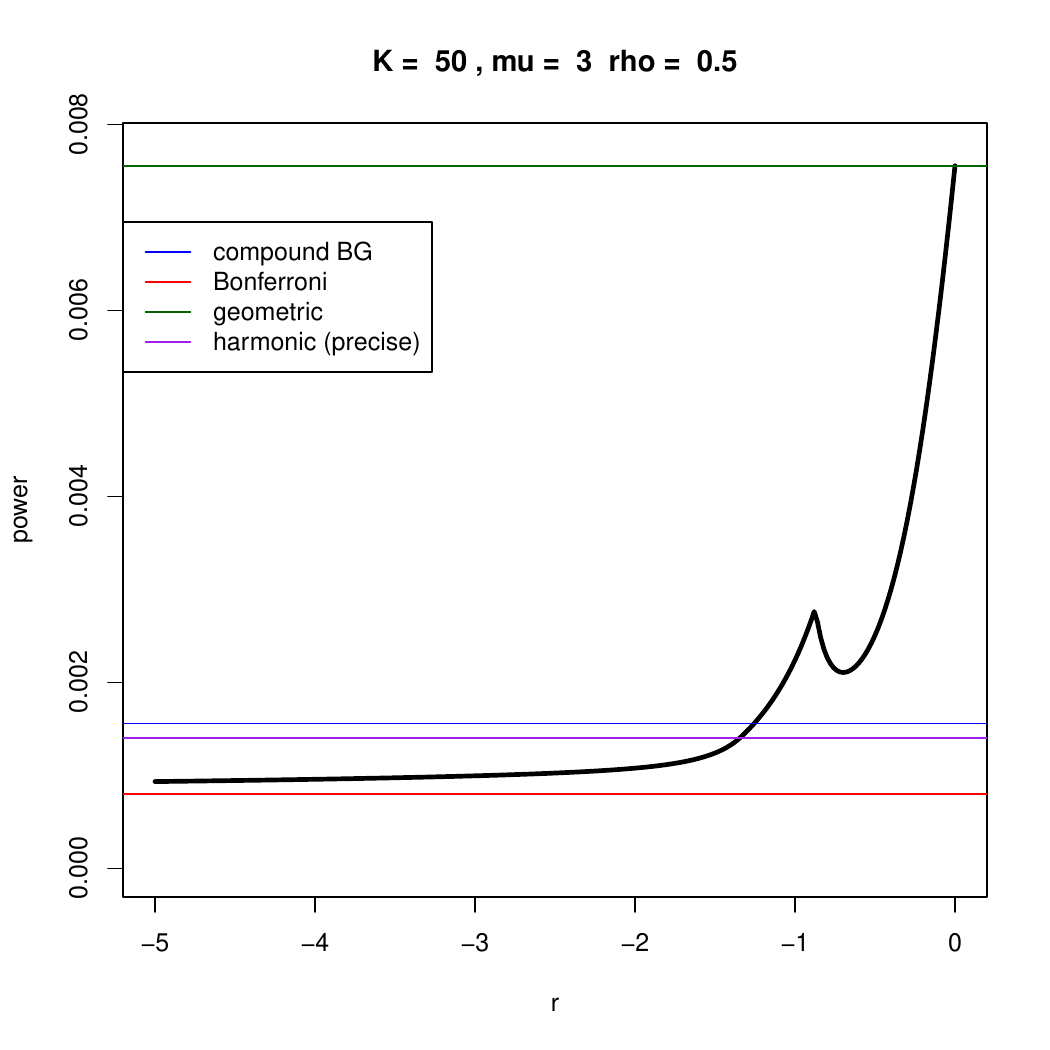} 
    \includegraphics[width=0.48\textwidth]{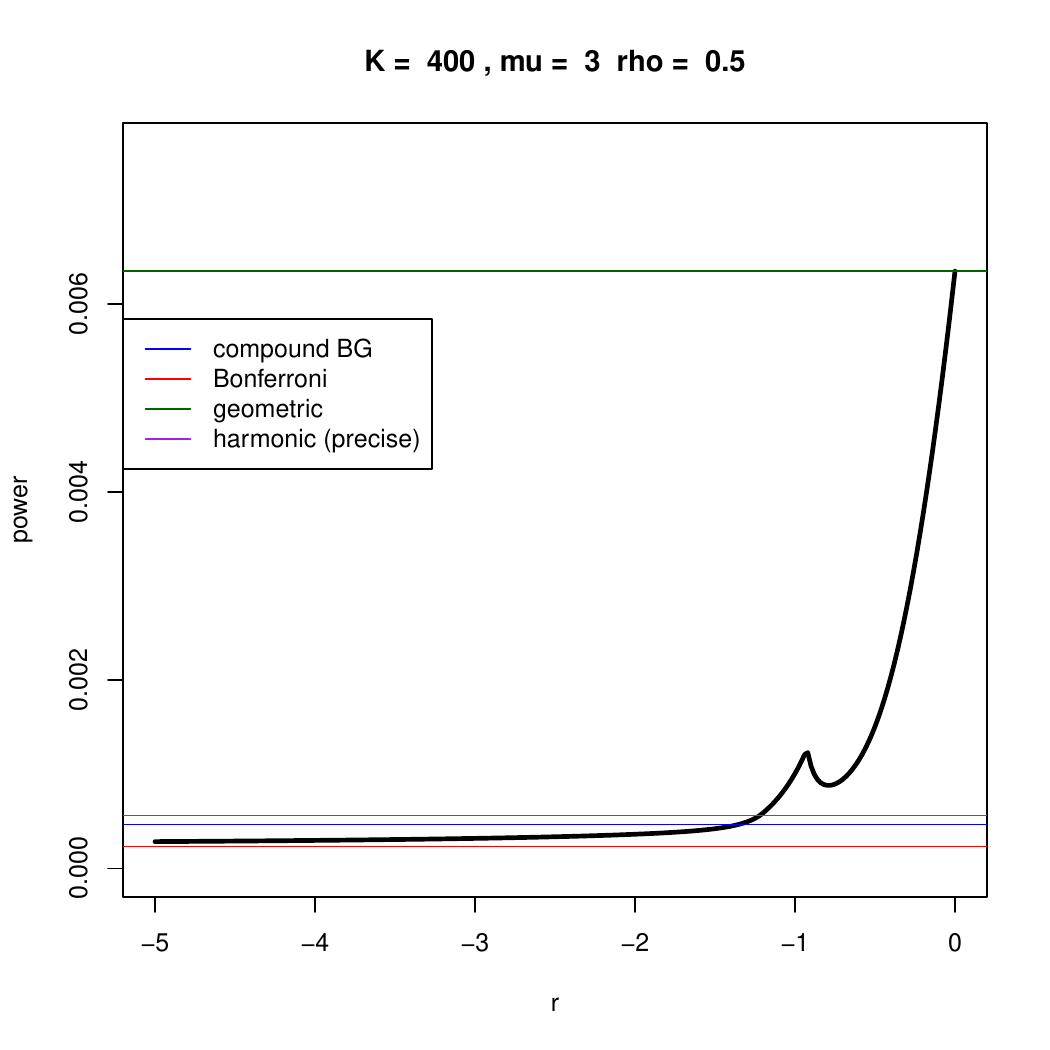}\\
    \includegraphics[width=0.48\textwidth]{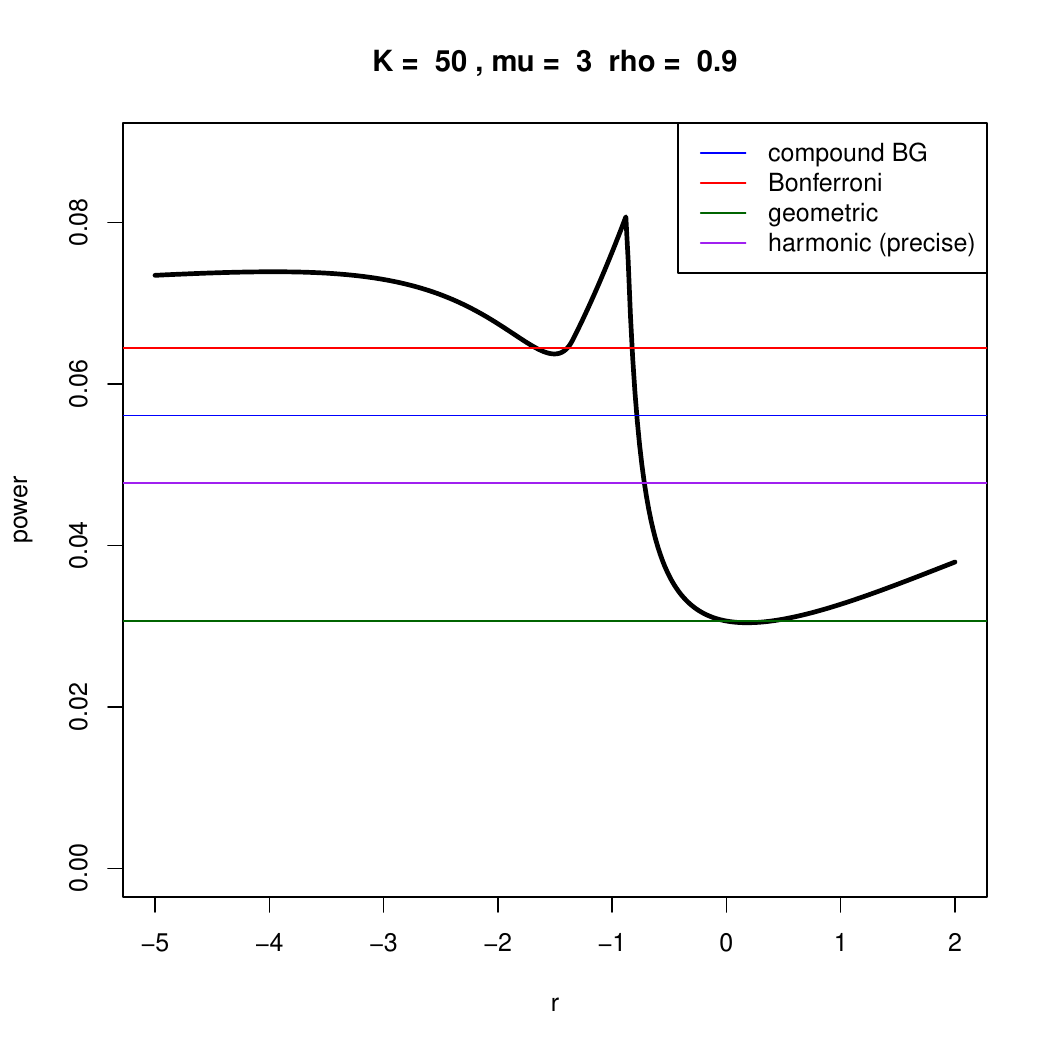}
    \includegraphics[width=0.48\textwidth]{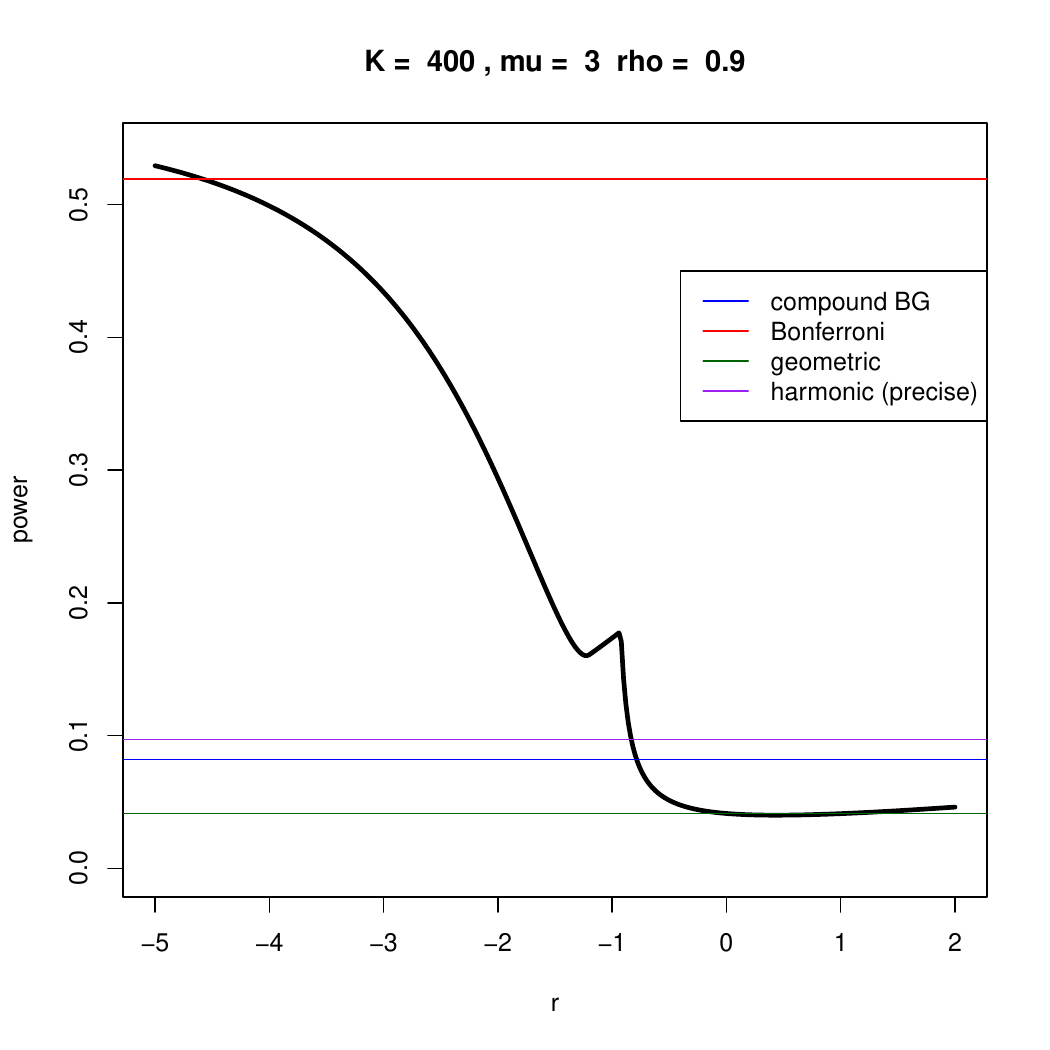}
  \end{center}
  \caption{Empirical averages of $P_{r,K}$  for $K\in \{50,400\}$}
  \label{fig:2}
\end{figure}

We conduct some simulations for the correlated standard z-tests in model \eqref{eq:defx} to compare different merging methods. 
In Figure \ref{fig:2}, we fix $\mu=3$ (other values of $\mu$ give qualitatively similar results),
and report the empirical performance for $K=50,400$ of merging methods with various values of $r$ as well as the compound Bonferroni-geometric method in \eqref{eq:BG}.
The curves in Figure \ref{fig:2} are based on an average of 1,000 replications of $P_{r,K}$. 
For the best of visibility, we plot the range $r\in[-5,0]$ for $\rho=0.1,0.5$ and $r\in[-5,2]$ for $\rho=0.9$. 
For values of $r$ close to $-1$, the asymptotically precise choice of the multiplier $a_{r,K}$ in Table~\ref{tab:main} explodes as $r$ approaches $-1$. 
To avoid this explosion for finite $K$, we replace $a_{r,K}$ by $e\ln K $ if $a_{r,K}>e\ln K$ and $r\ge \ln K/(1-\ln K)$,
and this gives valid merging functions  (cf.~\eqref{eq:min}; there is still a kink around $r=-1$ since these choices are not precise). 
We make the following observations from Figure~\ref{fig:2}.
\begin{enumerate} 
\item
  If the dependence is light or moderate ($\rho=0.1,0.5$), the Bonferroni method, as well as other methods based on $r<-1$,
  work quite well, and their power improves when $K$ increases from $50$ to $400$.
\item
  In the case of very strong dependence ($\rho=0.9$), the geometric averaging method and the arithmetic averaging method perform quite well.
  Notably, the performance of the Bonferroni method and other methods with $r<-1$ gets worse when $K$ increases from $50$ to $400$. 
\item
  The compound BG method generally performs quite well in all cases as it is often slightly worse (on average) than the best of the two base methods.
\end{enumerate}

\section{Directions of further research}
\label{sec:conclusion}

Perhaps the most important direction of further research is to find practically useful applications,
in multiple testing of a single hypothesis or testing multiple hypotheses,
for our methods of combining p-values.
The Bonferroni method of combining a set of p-values works very well
when experiments are almost independent,
while it produces unsatisfactory results under heavy dependence
(e.g., if all p-values are approximately equal).
Our methods are designed to work for  intermediate situations.
Promising results have been obtained using the harmonic mean \cite{Wilson:2019},
but other merging functions proposed in this paper also deserve careful experimental study.

Our emphasis has been on finding valid methods of combining p-values,
and we have just started exploration of their efficiency;
the results of Section~\ref{sec:efficiency} are asymptotic and crude.
A natural next step is to explore convergence of the combined p-values in distribution;
such results would be still asymptotic but more precise
(the price to pay will be the need to impose various conditions
on the distribution of the p-values).

The main application of results about efficiency of various combination methods
is choosing a suitable method.
Subsection \ref{subsec:search} is only a first step in this direction.
Much wider families of potential merging functions deserve to be explored.
In particular, finding the optimal value of $b$ in merging functions of the form \eqref{eq:adaptive}
is an interesting open problem.

This paper concentrates on the symmetric case,
where the merging function is a symmetric function of $p_1,\dots,p_K$.
However, the weighted case (as in Section~\ref{sec:weighted}) is important in many applications
because the quality of different p-values can be very different.
It would be of much interest to study efficient ways of assigning weights to the p-values
using prior or side information.

\subsection*{Acknowledgments}

We thank the referees and editors of the journal version of this paper for very constructive comments. 
We are grateful to Dave Cohen, Alessio Sancetta, Wouter Koolen, and Lutz Mattner for their advice.
Our special thanks go to Paul Embrechts, who, in addition to illuminating discussions,
has been instrumental in starting our collaboration.

For numeric calculations, we used Python and R.
The first author's work was partially supported by the Cyprus Research Promotion Foundation
and the European Union's Horizon 2020 Research and Innovation programme (671555),
and the second author's by the Natural Sciences and Engineering Research Council of Canada
(RGPIN-2018-03823, RGPAS-2018-522590).

\appendix

\numberwithin{equation}{section}
\setcounter{equation}{0}

\section{Proofs of main results}
\label{app:proofs-validity}

This appendix contains proofs of our main results. 
All these proofs are based on known results in robust risk aggregation.
For $a>0$, let $\mathcal{U}(a)$ be the set of all random variables distributed uniformly over the interval $[0,a]$, $a\ge0$;
we can regard $\mathcal{U}$ as an abbreviation for $\mathcal{U}(1)$.

\subsection*{Proof of  Lemma \ref{lem:RW-reduction}}  \label{pf:lem1}
  \textsc{Part ``if'' of (a):}
  Suppose $\qwc_\epsilon(F)\ge \epsilon$ for all $\epsilon\in (0,1)$.
  Consider arbitrary $U_1,\dots,U_K\in\mathcal{U}$.
  We have $q_\epsilon(F(U_1,\dots,U_K))\ge\epsilon$ for all $\epsilon\in(0,1)$.
  By the definition of left quantiles, $\Prob(F(U_1,\dots,U_K)<\epsilon)\le\epsilon$.
  It follows that, for all $\delta\in (0,1-\epsilon)$,
  \[
    \Prob(F(U_1,\dots,U_K)
    \le
    \epsilon)
    \le
    \Prob(F(U_1,\dots,U_K)<\epsilon+\delta)
    \le
    \epsilon+\delta,
  \]
  which implies
  \begin{equation*}
    \Prob(F(U_1,\dots,U_K)\le\epsilon)
    \le
    \epsilon,
  \end{equation*}
  since $\delta$ is arbitrary.
  Therefore, $F$ is a merging function.

  \textsc{Part ``only if'' of (a):}
  Suppose $F$ is a merging function.
  Let $U_1,\dots,U_K\in\mathcal{U}$ and $\epsilon\in(0,1)$.
  We have $\Prob(F(U_1,\dots,U_K)\le\epsilon)\le\epsilon$.
  By the definition of right quantiles, $\q^+_\epsilon (F(U_1,\dots,U_K))\ge\epsilon $.
  It follows that, for all $\delta\in(0,\epsilon)$,
  \[
    \q_\epsilon(F(U_1,\dots,U_K))
    \ge
    \q^+_{\epsilon-\delta}(F(U_1,\dots,U_K))
    \ge
    \epsilon-\delta,
  \]
  which implies $\q_\epsilon(F(U_1,\dots,U_K))\ge\epsilon$ since $\delta$ is arbitrary.

  \textsc{Part ``if'' of (b):}
  Suppose $\qwc_\epsilon(F)=\epsilon$ for all $\epsilon\in(0,1)$.
  By (a), $F$ is a merging function.
  For all $\epsilon,\delta \in (0,1)$,
  there exist $U_1,\dots,U_K\in\mathcal{U}$ such that $\q_\epsilon(F(U_1,\dots,U_K)) \in [\epsilon,\epsilon+\delta)$,
  which implies
  $\Prob(F(U_1,\dots,U_K)\le\epsilon+\delta)\ge \epsilon$.
  Since $\delta$ is arbitrary, we have
  \[
    \sup
    \left\{
      \Prob(F(U_1,\dots,U_K)\le\epsilon)\st U_1,\dots,U_K\in\mathcal{U}
    \right\}
    =
    \epsilon,
  \]
  and thus $F$ is precise.

  \textsc{Part ``only if'' of (b):}
  Suppose $F$ is a precise merging function.
  Since $F$ is a merging function, by (a) we have  $\qwc_\epsilon(M)\ge\epsilon$ for all $\epsilon\in(0,1)$.
  Suppose, for the purpose of contradiction, that $\qwc_\epsilon(M)>\epsilon$ for some $\epsilon\in(0,1)$.
  Then, there exists $\delta\in(0,1-\epsilon)$ such that $\q_\epsilon(F(U_1,\dots,U_K))>\epsilon+\delta$
  for all $U_1,\dots,U_K\in\mathcal{U}$.
  As a consequence, we have
  \[
    \Prob(F(U_1,\dots,U_K) \le \epsilon+\delta/2)
    \le
    \Prob(F(U_1,\dots,U_K) < \epsilon+\delta)
    \le
    \epsilon.
  \]
  Therefore,
  \[
    \sup
    \left\{
      \Prob(F(U_1,\dots,U_K) \le \epsilon+\delta/2) \st U_1,\dots,U_K\in\mathcal{U}
    \right\}
    <
    \epsilon+\delta/2,
  \]
  contradicting $F$ being precise.
 
In all proofs below, for statements that have a weighted version in Section~\ref{sec:weighted},
namely Theorem~\ref{thm:Bernard} and Propositions~\ref{prop:precise} and \ref{prop:super-precise},
we present a proof of the corresponding weighted version, which is stronger.

\subsection*{Proof of Theorem~\ref{thm:Bernard-w} (weighted version of Theorem~\ref{thm:Bernard})}

Without loss of generality we can, and will, assume that $\phi$ is strictly increasing.
Indeed, if $\phi$ is strictly decreasing, we can redefine $\phi:=-\phi$ and $\psi(u):=\psi(-u)$
and notice that the statement of the theorem for new $\phi$ and $\psi$
will imply the analogous statement for the original $\phi$ and $\psi$.

Define an accessory function $\Phi:(0,1)\to[-\infty,\infty]$
by $\Phi(\epsilon) = \frac{1}{\epsilon} \int_0^\epsilon \phi(u)\dd u$.
Fix $\epsilon\in(0,1)$.
Since $\phi$ is integrable, $\Phi(\epsilon)$ is finite.

Known results from the literature on robust risk aggregation can be applied
to random variables $X_k:=\phi(U_k)$, where $U_k\in\mathcal{U}$;
notice that the distribution function of $X_k$ is $\psi$:
\begin{equation*}
  \Prob(X_k\le x)
  =
  \Prob(\phi(U_k)\le x)
  =
  \Prob(U_k\le\psi(x))
  =
  \psi(x).
\end{equation*}
Theorem~4.6 of \cite{Bernard/etal:2014} gives the following relation:
\begin{equation}\label{eq:rw2}
  \qwc_{\epsilon}(M_{\phi,\w})
  =
  \inf
  \left\{
    \q_1
    \left(
      \psi
      \left(
        \sum_{k=1}^K w_k\phi(V_k)
      \right)
    \right)
    \st
    V_1,\dots,V_K\in\mathcal{U}(\epsilon)
  \right\}.
\end{equation}
Since
\[
  \q_1
  \left(
    w_1\phi(V_1)+\dots+w_K\phi(V_K)
  \right)
  \ge
  \Expect
  \left[
    w_1\phi(V_1)+\dots+w_K\phi(V_K)
  \right]
  =
  \Phi(\epsilon)
\]
for $V_1,\dots,V_K\in\mathcal{U}(\epsilon)$,
we have $\qwc_\epsilon(M_{\phi,\w}) \ge \psi(\Phi(\epsilon))$.

\subsection*{Proof of Proposition~\ref{prop:asymptotically-precise}}

The case $r=\infty$ is trivial, and we focus on $r\in (-1,\infty)$.
Let $\phi(u)=u^r$, which gives $\Phi(\epsilon)=\epsilon^r/(r+1)$,
in the notation of the previous proof.
Evaluating the term~\eqref{eq:term} in~\eqref{eq:general}, we obtain
\[
  \psi
  \left(
    \frac{1}{\epsilon}
    \int_0^{\epsilon}
      \phi(u)
    \dd u
  \right)
  =
  \begin{cases}
    \epsilon/e & \text{if $r=0$}\\
    (r+1)^{-1/r}\epsilon & \text{otherwise}.
  \end{cases}
\]
This shows that $(r+1)^{1/r} M_{r,K}$ is indeed a merging function.
  
Using Corollary~3.4 of~\cite{Embrechts/etal:2015}, we have
\[
  \lim_{K\to\infty}
  \frac{\phi(\qwc_\epsilon(M_{r,K}))}{\Phi(\epsilon)}
  =
  1,
\]
leading to 
\begin{equation}
  \lim_{K\to\infty} \qwc_\epsilon (M_{r,K}) = \epsilon(r+1)^{-1/r}. \label{eq:r1-2}
\end{equation}
It follows that, for $a<(r+1)^{1/r}$,
\[
  \lim_{K\to\infty}
  \qwc_{\epsilon}(a M_{r,K})
  <
  \epsilon
\]
and so, by Lemma \ref{lem:RW-reduction}, $a M_{r,K}$ is not a merging function for $K$ large enough.

\subsection*{Proof of Proposition~\ref{prop:precise-w}  (weighted version of Proposition~\ref{prop:precise})}

Let $M=(r+1)^{1/r}M_{r,\w}$.
Using \eqref{eq:rw2} and Theorem~\ref{thm:Bernard-w},
we have, for $\epsilon\in(0,1)$:
\begin{equation}\label{eq:rw3}
  \left(
    \qwc_\epsilon(M_{r,\w})
  \right)^r
  =
  \inf
  \left\{
    \q_1
    \left(
      \sum_{k=1}^K w_kV_k^r
    \right)
    \st
    V_1,\dots,V_K\in\mathcal{U}(\epsilon)
  \right\}
  \ge
  \frac{\epsilon^r}{1+r}
\end{equation}
if $r>0$,
\begin{equation}\label{eq:used-in-range}
  \left(
    \qwc_\epsilon(M_{r,\w})
  \right)^r
  =
  \sup
  \left\{
    \q_0^+
    \left(
      \sum_{k=1}^K w_kV_k^r
    \right)
    \st
    V_1,\dots,V_K\in\mathcal{U}(\epsilon)
  \right\}
  \le
  \frac{\epsilon^r}{1+r}
\end{equation}
if $r<0$, and
\begin{equation}\label{eq:rw4}
  \qwc_\epsilon(M_{r,\w}) 
  =
  \exp
  \left(
    \inf
    \left\{
      \q_1
      \left(
   \sum_{k=1}^K w_k \ln V_k
      \right)    \st
      V_1,\dots,V_K\in\mathcal{U}(\epsilon)
    \right\}
  \right)
  \ge
  \frac{\epsilon}{e}
\end{equation}
if $r=0$.
By Lemma \ref{lem:RW-reduction},
$M$ is a precise merging function if and only if the inequality in \eqref{eq:rw3}--\eqref{eq:rw4}
is an equality for all $\epsilon\in(0,1)$.

Fix $\epsilon\in(0,1)$ and $r\in(-1,\infty)$.
For $k=1,\dots,K$,
let $F_k$ be the distribution of $w_k V_k^r$ where $V_k\in\mathcal{U}(\epsilon)$.
Using the terminology of \cite{Wang/Wang:2016}, notice that the inequality in \eqref{eq:rw3}--\eqref{eq:rw4} is an equality
if and only if $(F_1,\dots,F_K)$ is jointly mixable 
due to a standard compactness argument (see \cite[Proposition 2.3]{Wang/Wang:2016}).
Therefore, we can first settle the cases $r=0$ and $r<0$,
as in these cases the supports of $F_1,\dots,F_K$ are unbounded on one side,
and $(F_1,\dots,F_K)$ is not jointly mixable (see \cite[Remark~2.2]{Wang/Wang:2016}).

Next assume $r>0$.
Since $F_1,\dots,F_K$ have monotone densities on their respective supports,
by Theorem~3.2 of \cite{Wang/Wang:2016},
$(F_1,\dots,F_K)$ is jointly mixable if and only if  the ``mean condition''
\[
  w\epsilon^r
  \le
  \frac{\epsilon^r}{1+r}
  \le
  \epsilon^r - w\epsilon^r
\]
is satisfied.
This is equivalent to $w\le\frac{1}{1+r}\le1-w$ and, therefore,
to the conjunction of $w\le1/2$ and $r\in[\frac{w}{1-w},\frac{1-w}{w}]$.
This completes the proof.

\subsection*{Proof of Proposition~\ref{prop:super-precise-w}  (weighted version of Proposition~\ref{prop:super-precise})}

Notice that, for each $k=1,\dots,K$,
the distribution of $w_k U_k^r$, where $U_k\in\mathcal{U}$, has a decreasing density on its support.
Therefore, we can apply Corollary~4.7 of \cite{Jakobsons/etal:2016},
which gives
\[
  \inf
  \left\{
    \q_\epsilon
    \left(
      w_1 U_1^r+\dots+w_K U_K^r
    \right)
    \st
    U_1,\dots,U_K \in\mathcal{U}
  \right\}
  =
  \max
  \left(
    w\epsilon^r,
    \frac{\epsilon^r}{1+r}
  \right).
\]
Simple algebra leads to
\[
  \qwc_\epsilon(M_{r,K})
  =
  \max
  \left(
    w, \frac{1}{1+r}
  \right)^{1/r}
  \epsilon,
\]
and by Lemma~\ref{lem:RW-reduction}, $M$ is a precise merging function.

\subsection*{Proof of Proposition~\ref{prop:geometric}}

First, we note that as $\lim_{r\to 0} (r+1)^{1/r}=e$, 
we know $a^G_K\le e$  from Proposition \ref{prop:asymptotically-precise}. 
Moreover, by letting $r\to 0$ in \eqref{eq:r1-2}, 
we know that $a^G_K\to e$.

Our next goal is to obtain the precise value of $\qwc_{\epsilon}(M_{0,K})$.
Set
\begin{align*}
  b_K
  &:=
  \sup
  \left\{
    \q_0^+
    \left(
      -(\ln U_1+\dots+\ln U_K)
    \right)
    \st
    U_1,\dots,U_K\in\mathcal{U}
  \right\}\\
  &=
  \sup
  \left\{
    \q_0^+
    \left(
      -(\ln V_1+\dots+\ln V_K) + K\ln\epsilon
    \right)
    \st
    V_1,\dots,V_K\in \mathcal{U}(\epsilon)
  \right\}.
\end{align*}
It is easy to see that
\begin{align*}
  \qwc_{\epsilon}(M_{0,K})
  &=
  \exp
  \left(
    \inf
    \left\{
      \q_1
      \left(
        \frac{\ln V_1+\dots+\ln V_K}{K}
      \right)
      \st
      V_1,\dots,V_K\in\mathcal{U}(\epsilon)
    \right\}
  \right)\\
  &=
  \exp
  \left(
    -\sup
    \left\{
      \q_0^+
      \left(
        -\frac{\ln V_1+\dots+\ln V_K}{K}
      \right)
      \st
      V_1,\dots,V_K\in\mathcal{U}(\epsilon)
    \right\}
  \right)\\
  &=
  \exp
  \left(
    -b_K/K+\ln\epsilon
  \right)\\
  &=
  \epsilon
  \exp(-b_K/K).
\end{align*}
It is clear that $e^{b_K/K}M_{0,K}$ is a precise merging function, and
$b_K/K\to 1$.

Since $-\ln U$ has the standard exponential distribution for $U\in\mathcal{U}$
and, therefore, a decreasing density on $\R$,
we can apply Theorem~3.2 of \cite{Bernard/etal:2014}
(essentially Theorem~3.5 of \cite{Wang/Wang:2011})
to arrive at
\[
  b_K
  =
  -(K-1)\ln(1-(K-1)c_K)-\ln c_K,
\]
where $c_K$ is the unique solution to \eqref{eq:ga} (see \cite[Corollary~4.1]{Wang/Wang:2011}).
Using \eqref{eq:ga}, we can write
\[
  b_K/K
  =
  - \ln c_K - (K-1)(1-K c_K).
\]
Using $a^G_K=e^{b_K/K}$ one obtains the desired result.

\subsection*{Proof of Theorem~\ref{thm:Embrechts}}

We will apply Theorem~4.2 of \cite{Embrechts/Puccetti:2006} in our situation
where the function $\phi$ (and, therefore, $\psi$ as well) in~\eqref{eq:mean}
is decreasing.
Letting $X_k:=\phi(p_k)$ and using the notation $m_+$
(used in Theorem~4.2 of \cite{Embrechts/Puccetti:2006}),
we have, by the definition of $m_+$,
\begin{align*}
  \Prob
  \left(
    \sum_{k=1}^K X_k
    <
    s
  \right)
  &\ge
  m_+(s),\notag\\
  \Prob
  \left(
    \frac1K \sum_{k=1}^K \phi(p_k)
    <
    s/K
  \right)
  &\ge
  m_+(s),\notag\\
  \Prob
  \left(
    M_{\phi,K}(p_1,\dots,p_K)
    >
    \psi(s/K)
  \right)
  &\ge
  m_+(s),\notag\\
  \Prob
  \left(
    M_{\phi,K}(p_1,\dots,p_K)
    \le
    \psi(s/K)
  \right)
  &\le
  1-m_+(s).
\end{align*}
The lower bound on $m_+(s)$ given in Theorem~4.2 of \cite{Embrechts/Puccetti:2006}
involves $1-F(x)$, where $F$ is the common distribution function of $X_k$,
and in our current context we have:
\begin{equation*}
  1-F(x)
  =
  \Prob(X_k>x)
  =
  \Prob(\phi(p_k)>x)
  =
  \Prob(p_k<\psi(x))
  =
  \psi(x).
\end{equation*}
The last inequality and chain of equalities
in combination with Theorem~4.2 of \cite{Embrechts/Puccetti:2006}
give
\[
  \Prob
  \left(
    M_{\phi,K}(p_1,\dots,p_K)
    \le
    \psi(s/K)
  \right)
  \le
  K
  \inf_{r\in[0,s/K)}
  \frac{\int_r^{s-(K-1)r}\psi(x)\dd x}{s-Kr}.
\]
Setting $\epsilon:=\psi(s/K)\in[\psi(\infty),\psi(0)]$
(so that it is essential that $\psi(\infty)=0$),
we obtain, using $s=K\phi(\epsilon)$,
\[
  \Prob
  \left(
    M_{\phi,K}(p_1,\dots,p_K)
    \le
    \epsilon
  \right)
  \le
  K
  \inf_{r\in[0,\phi(\epsilon))}
  \frac{\int_r^{K\phi(\epsilon)-(K-1)r}\psi(x)\dd x}{(\phi(\epsilon)-r)K}.
\]
Setting $t:=\phi(\epsilon)-r$ and renaming $x$ to $u$,
this can be rewritten as~\eqref{eq:Embrechts}.
 
\subsection*{Proof of Proposition~\ref{prop:non-integrable-optimality}}

The case $r=-\infty$ is trivial, and we focus on $r\in (-\infty,-1)$. 
By Theorem~\ref{thm:Embrechts} applied to $\phi(u):=u^r$, $r<-1$, we have 
\begin{equation*}
  \Prob
  \left(
    M_{\phi,K}(p_1,\dots,p_K)
    \le
    \epsilon
  \right)
  \le
  \frac{\int_{0}^{K\phi(\epsilon)}\psi(u)\dd u}{\phi(\epsilon)}
  =
  \frac{r}{r+1}
  K^{1+1/r}
  \epsilon. 
\end{equation*}
This shows that
$
  \frac{r}{r+1}
  K^{1+1/r}
  M_{r,K}
$
is indeed a merging function.

We next show a simple property of a precise merging function via general averaging.
Define the following constant:
\[
  b_{r,K}
  :=
  \left(
    \frac{1}{K}
    \sup\{\q_0^+(U_1^r+\dots+U_K^r)~ \big |~ U_1,\dots,U_K\in\mathcal{U}\}
  \right)^{-1/r}.
\]
It is clear that $b_{r,K}\ge 1$  for $r<0$.

\begin{lemma}\label{lem:lem13}
  For $r<0$, the function $b_{r,K} M_{r,K}$ is a precise merging function.
\end{lemma}
\begin{proof}
  By straightforward algebra and Theorem~4.6 of \cite{Bernard/etal:2014},
  \begin{align*}
    \qwc_\epsilon(b_{r,K} M_{r,K})
    &=
    b_{r,K}
    \inf
    \left\{
      \q_1(M_{r,K}(V_1,\dots,V_K)) \st V_1,\dots,V_K\in\mathcal{U}(\epsilon)
    \right\}\\
    &=
    b_{r,K}
    \inf
    \left\{
      \epsilon \q_1(M_{r,K}(U_1,\dots,U_K)) \st U_1,\dots,U_K\in\mathcal{U}
    \right\}\\
    &=
    b_{r,K} \epsilon
    \left(
      \frac1K
      \sup
      \left\{
        \q_0^+(U_1^r+\dots+U_K^r) \st U_1,\dots,U_K\in\mathcal{U}
      \right\}
    \right)^{1/r}\\
    &=
    \epsilon.
  \end{align*}
  By Lemma~\ref{lem:RW-reduction}, $b_{r,K} M_{r,K}$ is a precise merging function.
\end{proof}

To construct precise merging functions, it remains to find values of $b_{r,K}$.
Unfortunately, for $r<0$ no analytical formula for $b_{r,K}$ is available.
There is an asymptotic result available in \cite{Bignozzi/etal:2016},
which leads to the following proposition.

\begin{proposition}\label{prop:rw6}
  For $r\in (-\infty,-1)$,
  \[
    \lim_{K\to\infty}
    \frac{b_{r,K}}{K^{1+1/r}}
    =
    \frac{r}{r+1}.
  \]
\end{proposition}

\begin{proof}
  The quantity $\overline{\Delta}^{\mathcal{F}_d}$ in \cite{Bignozzi/etal:2016}, defined as
  \[
    \overline{\Delta}^{\mathcal{F}_d}
    :=
    \lim_{\alpha\to1}
    \frac
    {
      \sup
      \left\{
        \q_\alpha(U_1^r+\dots+U_K^{r}) \st U_1,\dots,U_K\in\mathcal{U}(\alpha)
      \right\}
    }
    {K(1-\alpha)^r},
  \]
  satisfies
  \begin{equation*}
    \overline{\Delta}^{\mathcal{F}_d}
    =
    \frac{1}{K}
    \sup
    \left\{
      \q_0^+(U_1^r+\dots+U_K^r) \st U_1,\dots,U_K\in\mathcal{U}
    \right\}
    =
    b_{r,K}^{-r}.
  \end{equation*}
  Using Proposition~3.5 of \cite{Bignozzi/etal:2016}, we have, for $r<-1$,
  by substituting $\beta:=-1/r$ in~(3.25) of \cite{Bignozzi/etal:2016}
  and $\overline{\Delta}^{\mathcal{F}_d}=b_{r,K}^{-r}$,
  \[
    \lim_{K\to\infty}
    \frac{b_{r,K}^{-r}}{K^{-r-1}}
    =
    \left(
      \frac{r}{r+1}
    \right)^{-r},
  \]
  and this gives the desired result.
\end{proof}
The claim on the asymptotic precision in Proposition~\ref{prop:non-integrable-optimality} immediately follows from Lemma \ref{lem:lem13} and Proposition~\ref{prop:rw6}.

\subsection*{Proof of Proposition~\ref{prop:harmonic}}

Let us first find the smallest value of the coefficient
$\frac{r}{r+1}K^{1+1/r}$ in  Proposition~\ref{prop:non-integrable-optimality}.
Setting the derivative in $r$ of the logarithm of this coefficient to 0,
we obtain a linear equation whose solution is
\begin{equation}\label{eq:min}
  r
  =
  \frac{\ln K}{1-\ln K}.
\end{equation}
Plugging this into the coefficient gives $e\ln K$.
Notice that $r$ defined by~\eqref{eq:min} satisfies $r<-1$
and apply the inequality $M_{r,K}\le M_{-1,K}$, a special case of  \eqref{eq:r1-1}.
Hence, $e\ln K M_{-1,K}$ is a merging function,
and $a^H_K\le e\ln K$ follows from the fact that $a^H_K M_{-1,K}$ is a precise merging function, which we prove below. 
   
Using the notation in Lemma \ref{lem:lem13}, we need to show $a^H_K=b_{-1,K}$, and we shall use Corollary~3.7 of \cite{Wang/etal:2013}. 
Write
\[
  H(t)
  :=
  \frac{K-1}{1-(K-1)t} + \frac1t
  =
  \frac{1}{t(1-(K-1)t)},
  \quad
  t\in [0,1/K].
\]
By Corollary~3.7 of \cite{Wang/etal:2013}, we have
\begin{align*}
 b_{-1,K}
  &=
  \frac{1}{K}
  \sup
  \left\{
    q_0^+(U_1^{-1}+\dots+U_K^{-1}) \st U_1,\dots,U_K\in \mathcal{U}
  \right\}\\
  &=
  \frac{1}{K} H(x_K)
  =
  \frac{1}{ Kx_K(1-(K-1)x_K)}
\end{align*}
where $x_K$ solves the equation
\[
  \int_x^{1/K} H(t)\dd t
  =
  \left(
    \frac1K - x
  \right)
  H(x),
  \quad
  x \in [0,1/K).
\]
Plugging in the expression for $H$ and rearranging the above equation, we obtain
\begin{align}
  \frac{1-Kx}{Kx(1-(K-1)x)}
  &=
  \int_x^{1/K}
  \left(
    \frac{K-1}{1-(K-1)t} + \frac 1t
  \right)
  \dd t \notag\\
  &=
  \int_{(K-1)x}^{(K-1)/K}
  \frac{1}{1-y}\dd y + \int_x^{1/K}\frac1t
  \dd t \notag\\
  &=
  \ln(1-(K-1)x)-\ln x.
  \label{eq:harmo2}
\end{align}
The uniqueness of the solution $x_K$ to \eqref{eq:harmo2} can be easily checked,
and it is a special case of Lemma~3.1 of \cite{Jakobsons/etal:2016}.
Writing $y=\frac{1-Kx}{x}>0$, \eqref{eq:harmo2} reads as
$
  \frac{y}{(y+1)K/(y+K)}  = \ln (y+1)
$.
Rearranging the terms gives
\begin{equation}\label{eq:harmo3}
  y^2
  =
  K( (y+1)\ln (y+1)-y),
\end{equation}
which admits a unique solution, $y_K=\frac{1-Kx_K}{x_K}$.
Therefore,
\[
  b_{-1,K}
  =
  \frac{1}{K x_K (1-(K-1)x_K)}
  =
  \frac{(y_K+K)^2}{(y_K+1)K},
\]
and hence $a^H_K  = b_{-1,K}$.

Next we analyze the asymptotic behaviour of $a_{ K}^H$ as $K\to\infty$.
Using $\ln(y+1)\ge y-y^2/2$ for $y \ge 0$,
we can see that \eqref{eq:harmo3} implies the inequality
\[
  y^2
  \ge
  \frac{K}{2} y^2 - \frac{K}{2} y^3,
\]
which leads to $2 \ge K(1-y)$.
Hence, we have  $\liminf_{K\to\infty} y_K \ge 1$.

Notice that $(y+1)\ln(y+1)-y$ is a strictly increasing function of $y\in(0,\infty)$.
Using $\liminf_{K\to \infty} y_K \ge 1$, we obtain that
\begin{equation*}
  \liminf_{K\to\infty} y^2_K
  \ge
  K(2\ln 2-1).
\end{equation*}
Therefore, $\lim_{K\to\infty}y_K=\infty$.
Applying logarithms to both sides of \eqref{eq:harmo3} and taking a limit in their ratio, we obtain
\begin{equation*}
  1
  =
  \lim_{K\to\infty}
  \frac
  {
    2\ln y_K
  }
  {
    \ln K + \ln (y_K+1)
    +
    \ln
    \left(
      \ln (y_K+1)-\frac{y_K}{y_K+1}
    \right)
  }
  =
  \lim_{K\to\infty}
  \frac{2\ln y_K}{\ln K + \ln y_K},
\end{equation*}
and hence $\ln y_K/\ln K\to 1$ as $K\to\infty$.
Using \eqref{eq:harmo3} again, we have
\begin{equation*}
  1
  =
  \lim_{K\to\infty}
  \frac{y_K^2}{K((y_K+1)\ln(y_K+1)-y_K)}
  =
  \lim_{K\to\infty}
  \frac{y_K^2}{ K( y_K\ln y_K)}
  =
  \lim_{K\to\infty}
  \frac{y_K}{K\ln K}.
\end{equation*}
Therefore, we have
\[
  \lim_{K\to\infty}
  \frac{a^H_{K}}{\ln K}
  =
  \lim_{K\to\infty}
  \frac{(y_K+K)^2}{(y_K+1)K \ln K}
  =
  \lim_{K\to\infty}
  \frac{(K\ln K)^2}{(K\ln K)K \ln K}
  =
  1.
\]
This completes the proof.

\subsection*{Proof of Proposition~\ref{prop:homo}}

Let us check that for $F:=a M_{r,K}$ the following statements are equivalent:
\begin{itemize}
\item[(a)]
  $F$ is a merging function;
\item[(b)]
  \eqref{eq:merging} holds for some $\epsilon\in(0,1)$;
\item[(c)]
  $\qwc_{\epsilon}(F)\ge\epsilon$ for some $\epsilon\in(0,1)$.
\end{itemize}

The implication $\text{(a)}\Rightarrow\text{(b)}$ holds by definition.

To check $\text{(b)}\Rightarrow\text{(c)}$,
let us assume (b).
Since $\Prob(X\le \epsilon)\le\epsilon$ implies $q_\epsilon^+(X)\ge\epsilon$ for any random variable~$X$,
\[
  \inf
  \left\{
    q^+_\epsilon(F(U_1,\dots,U_K)) \st U_1,\dots,U_K\in\mathcal{U}
  \right\}
  \ge
  \epsilon.
\]
Using Lemma~4.5 of \cite{Bernard/etal:2014},
\[
  \qwc_{\epsilon}(F)
  =
  \inf
  \left\{
    q^+_\epsilon(F(U_1,\dots,U_K)) \st U_1,\dots,U_K\in\mathcal{U}
  \right\},
\]
and hence (c) follows.

It remains to check $\text{(c)}\Rightarrow\text{(a)}$.
For any $ \epsilon\in(0,1)$, by straightforward algebra and Theorem~4.6 of \cite{Bernard/etal:2014},
\begin{align*}
  \qwc_{\epsilon}(F)
  &=
  \inf
  \left\{
    q_1(F(V_1,\dots,V_K)) \st V_1,\dots,V_K\in \mathcal{U}(\epsilon)
  \right\}\\
  &=
  \epsilon
  \inf
  \left\{
    q_1(F(U_1,\dots,U_K)) \st U_1,\dots,U_K\in\mathcal{U}
  \right\}.
\end{align*}
Therefore, to check $\qwc_{\epsilon}(F)\ge\epsilon$ for all $\epsilon\in(0,1)$,
one only needs to check the inequality for one $\epsilon\in(0,1)$.
By Lemma \ref{lem:RW-reduction}, $F$ is a merging function.

This completes the proof of the first part of Proposition~\ref{prop:homo},
and the second part can be proved similarly.

\subsection*{Proof of Proposition~\ref{prop:IID}}

First we consider the integrable case $r>-\Pi$.
By the Marcinkiewicz strong law of large numbers \cite[Theorem~6.7.1]{Gut:2013}
(and Kolmogorov's strong law of large numbers, which is its special case),
we have
\[
  \frac1K
  \sum_{k=1}^K
  p_k^r
  \to
  \Expect(p_1^r)
  \quad
  \text{a.s.}
\]
This implies parts (i) and (ii) of the case $\Pi>1$ and part (i) of the case $\Pi<1$.
Since $\Prob(p_1>0)=1$ unless $\Pi=0$, this also gives the ``finite positive'' entries
in Figure~\ref{fig:summary}.

Now let us consider the non-integrable case $r<-\Pi$.
In this case $-\Pi/r\in[0,1)$ and we can apply Marcinkiewicz's strong law,
assuming $\Pi>0$ (the simple case $\Pi=0$ should be considered separately).
We have $\Expect((p_1^r)^{-\Pi/r-\epsilon})<\infty$ for any $\epsilon>0$,
and so, by Marcinkiewicz's law (applicable when $\epsilon<-\Pi/r$),
\[
  \frac{\sum_{k=1}^K p_k^r}{K^{-r/\Pi+\epsilon}}
  \to
  0
  \quad
  \text{a.s.}
\]
for any $\epsilon>0$.
This implies
\begin{equation}\label{eq:side-1}
  \frac{\left(\sum_{k=1}^K p_k^r / K\right)^{1/r}}{K^{-1/\Pi-1/r-\epsilon}}
  \to
  \infty
  \quad
  \text{a.s.}
\end{equation}
for any $\epsilon>0$.
On the other hand,
we have $\Expect((p_1^r)^{-\Pi/r+\epsilon})=\infty$ for any $\epsilon>0$,
and so, by the other part of Marcinkiewicz's law (applicable when $\epsilon<1+\Pi/r$),
\[
  \frac{\sum_{k=1}^K p_k^r}{K^{-r/\Pi-\epsilon}}
  \to
  0
  \quad
  \text{a.s.}
\]
fails for any $\epsilon>0$.
Therefore,
\[
  \frac{\left(\sum_{k=1}^K p_k^r / K\right)^{1/r}}{K^{-1/\Pi-1/r+\epsilon}}
  \to
  \infty
  \quad
  \text{a.s.}
\]
fails for any $\epsilon>0$.
This implies that
\begin{equation}\label{eq:side-2}
  \liminf_{K\to\infty}
  \frac{\left(\sum_{k=1}^K p_k^r / K\right)^{1/r}}{K^{-1/\Pi-1/r+\epsilon}}
  =
  0
\end{equation}
with a positive probability for any $\epsilon>0$.
Kolmogorov's zero-one law \cite[Theorem 1.5.1]{Gut:2013}
allows us to replace ``with a positive probability'' by ``a.s.''
Combining \eqref{eq:side-1} and \eqref{eq:side-2} gives
part (iii) of the case $\Pi>1$ and
parts (ii) and (iii) of the case $\Pi<1$.

It remains to consider the case $r=-\infty$ (Bonferroni).
First suppose $\Pi>1$.
{For any $b\in (0, 1-1/\Pi)$ and $B>0$,  take $c = (1/(1-b),\Pi)$, and note that  $ \Expect(p_1^{-c})<\infty$.
We have
\begin{align*}
  \Prob(P_{r,K} \le B K^b)
  &=
  \Prob(\exists k\in\{1,\dots,K\}: p_k \le B  K^{b-1})
  \\& \le K \Prob(p_1 \le B K^{b-1}) 
  \\&= K \Prob(p_1^{-c} \ge (B K^{b-1} )^{-c})
  \\& \le  K   \Expect(p_1^{-c})  / (B K^{b-1})^{-c} 
   =
  \Expect(p_1^{-c}) B^{c}  K^{1-(1-b)c}\to 0.
\end{align*}
This gives $P_{r,K}/K^b\to\infty$ in probability and, therefore,
part (iv) of the case $\Pi>1$.}
Finally, suppose $\Pi<1$. 
Note that
for any $r<-1$, 
\[
  P_{-\infty, K} = K \min_{k\in \{1,\dots,K\}} p_k \le K M_{r,K} (p_1,\dots,p_K) =   \frac{r+1}{r }K^{-1/r} P_{r,K}.
\]
For $b>1-1/\Pi$, we can take $r<-1$ small enough such that 
$ b+1/r>1-1/\Pi$, and obtain by using (iii), 
\[
  P_{-\infty, K}/K^b \le   \frac{r+1}{r } P_{r,K}/ K^{b+1/r} \to 0 ~~~\text{a.s.}
\]
 
\subsection*{Proof of Proposition~\ref{prop:BA-BG}}

We first analyze the compound BA method. 
To show that the family $F_K^{\BA}$ is asymptotically precise,
it suffices to show that for some $\epsilon\in(0,1)$,
\begin{equation}\label{eq:adaptive-3}
  \limsup_{K\to\infty}
  \left(
    \sup
    \left\{
      \Prob(F_K^{\BA}(U_1,\dots,U_K)\le\epsilon)
      \st
      U_1,\dots,U_K\in\mathcal{U}
    \right\}
  \right)
  =
  \epsilon.
\end{equation} 
To show \eqref{eq:adaptive-3}, fix $\epsilon \in (0,1)$ and an even number $K$.
Set $\delta:=\epsilon/(2K)$.
Let $A_1,\dots,A_K,B$ be disjoint events with $\Prob(A_j)= \delta$, $j=1,\dots,K$,
and $\Prob(B)= (K-2) \delta$.
Let $V, V_1,\dots,V_K, W_1,\dots,W_K $ be independent random variables,
between themselves and of $A_1,\dots,A_K,B$,
such that, for $j=1,\dots,K$,
$V_j$ is uniformly distributed on $[0,\delta]$,
$V$ is uniformly distributed on $[\delta, (K-1)\delta]$,
and $W_j$ is uniformly distributed on $[(K-1)\delta, 1]$.

Define the random variables
\[
  p_j
  :=
  \III_{B} V + \III_{A_j}  V_j + \III_{(B\cup A_j)^\complement} W_j
\]
for an odd number $j\in\{1,\dots,K\}$,
and
\[
  p_j
  :=
  \III_{B} (K\delta -V) +\III_{A_j}  V_j + \III_{(B\cup A_j)^\complement} W_j
\]
for an even number $j\in\{1,\dots,K\}$.
We can easily check that $p_1,\dots,p_K$ are uniformly distributed on $[0,1]$.
Note that $\frac2K \sum_{k=1}^K p_k \III_B = K\delta \III_B$,
and, by writing $A:=\bigcup_{j=1}^n A_K$,
$\min(p_1,\dots,p_K)\III_A \le \delta \III_A$.
Therefore,
\begin{align*}
 \Prob\left(F^{\BA}_K(p_1,\dots,p_K) \le \epsilon \right)
  &=
  \Prob
  \left(
    \min ( K\min(p_1,\dots,p_K), 2 \bar p) \le K\delta
  \right)\\
  &\ge
  \Prob(A \cup B)
  =
  (2K-2) \delta
  =
  \epsilon- \epsilon/K .
\end{align*} 
Hence, \eqref{eq:adaptive-3} holds, and the family $F^{\BA}_K$ is asymptotically precise.

The case of the compound BG method is similar,
although an explicit construction is more complicated.
We fix $\epsilon\in(0,1/2)$ and let $\delta := \epsilon/K$.
Note that, for a uniform random variable $U$ on $[0,1]$, we have
\[
  \Expect
  \left(
    \ln U \givn U \in [0,\epsilon]
  \right)
  =
  \ln\epsilon - 1.
\]
For $K \ge 3$, $\ln\delta < \ln\epsilon-1$,
and there exists a number $t(K)>0$ such that
\[
  \Expect
  \left(
    \ln U \givn U \in [\delta,\epsilon-t(K)]
  \right)
  =
  \ln\epsilon - 1.
\]
Clearly, $t(K)\to 0$ as $K\to\infty$.
Next, take an arbitrary $\eta>0$.
By \cite[Corollary 3.4]{Embrechts/etal:2015},
there exists $K_0\in\{1,2,\dots\}$ such that for $K>K_0$,
there exist uniform random variables $U_1,\dots,U_K$ on $[\delta,\epsilon-t(K)]$ satisfying
\[
  q_{1-\eta} (\ln U_1+\dots+\ln U_K)
  <
  \Expect(\ln U_1+\dots+\ln U_K)
  =
  K(\ln\epsilon - 1).
\]
In other words,
\[
  \Prob
  (e \tilde U \le \epsilon ) \ge 1-\eta,
\]
where $\tilde U$ is the geometric mean of $U_1,\dots,U_K$.

Let $A_1,\dots,A_K,B$ be disjoint events
with $\Prob(A_j)= \delta$, $j=1,\dots,K$, and $\Prob(B)= \epsilon - t (K)$.
Further, let $V_1,\dots,V_K, W_1,\dots,W_K$ be random variables,
independent between themselves and of $ U_1,\dots,U_K $ and $A_1,\dots,A_K,B$,
and such that, for $j=1,\dots,K$,
$V_j$ is uniformly distributed on $[0,\delta]$
and $W_j$ is uniformly distributed on $[\epsilon-t(K),1]$.

Define the following random variables, for $j=1,\dots,K$:
\[
  p_j
  :=
  \III_{B} U_j +\III_{A_j}  V_j + \III_{(B\cup A_j)^\complement}  W_j.
\]
We can easily check that $p_1,\dots,p_K$ are uniformly distributed on $[0,1]$. 
Note that
\[
  \Prob(e \tilde p \le \epsilon \givn B)
  =
  \Prob(e \tilde U \le \epsilon \givn B)
  \ge
  1-\eta.
\]
By setting $A:=\bigcup_{j=1}^n A_K$, we have
$\min(p_1,\dots,p_K)\III_A \le \delta \III_A$.
Therefore,
\begin{align*}
  \Prob
  \left(
    F^{\BG}_K(p_1,\dots,p_K) \le  2 \epsilon
  \right)
  &= 
  \Prob
  \left(
    \min(K\min(p_1,\dots,p_K), e\tilde p) \le \epsilon
  \right)\\
  &\ge
  \Prob(A ) + \Prob(B) (1-\eta)
  =
  \epsilon + (\epsilon - t(K))(1-\eta)  .
\end{align*} 
Since $\eta$ is arbitrary and $t(K)\to 0$ as $K\to\infty$,
we know that
\[
  \lim _{K\to \infty}
  \left(
    \sup
    \left\{
      \Prob(F_K^{\BG}(U_1,\dots,U_K)\le 2 \epsilon)
      \st
      U_1,\dots,U_K\in\mathcal{U}
    \right\}
  \right)
  =
  2\epsilon,
\] 
and so the family $F^{\BG}_K$ is asymptotically precise.

\begin{lemma}\label{lem:gaussian}
  Let $Z$ be a standard normal random variable and $z\in\R$.
  For all $\sigma>0$,
  $\E[(\Phi(\sigma Z + z))^r] < \infty$ if and only if $r>-1/\sigma^2$ or both $r=-1/\sigma^2$ and $z>0$.
\end{lemma}
\begin{proof}
  It suffices to investigate the case $r<0$.
  In this case, only small values of $\Phi(\sigma Z + z)$ matter.
  To analyze whether
  \begin{equation}\label{eq:gaussianint}
    \E[(\Phi(\sigma Z + z))^r]
    =
    \int_{-\infty}^{\infty} (\Phi(\sigma x + z))^r \frac{1}{\sqrt{2\pi}} e^{-x^2/2}\dd x
  \end{equation}
  is finite, it suffices to look at the limiting behavior of $\Phi(\sigma x  + z) $
  as $x \to -\infty$.
  L'H\^opital's rule gives
  $\Phi(y) \sim - \frac{1}{y} \frac{1}{\sqrt{2\pi}} e^{-y^2/2}$ as $y \to -\infty$.
  Hence, we have, for some constant $C>0$, as $x\to-\infty$,
  \begin{align*}
    (\Phi(\sigma x + z))^r e^{-x^2/2}
    &\sim
    \left(-\frac{1}{\sigma x + z} \frac{1}{\sqrt{2\pi}} \right)^r e^{-r(\sigma x + z)^2/2 - x^2/2}\\
    &\sim
    C x^{-r} e^{-(r\sigma^2+1)x^2/2 + r\sigma z x}.
  \end{align*} 
  Therefore, if $r\sigma^2+1>0$, the integral in \eqref{eq:gaussianint} is finite,
  and if $r\sigma^2+1<0$, the integral in \eqref{eq:gaussianint} is infinite.
  If $r\sigma^2+1=0$, then $z\le 0$ leads to an infinite integral in \eqref{eq:gaussianint},
  and $z>0$ leads to  a finite integral in \eqref{eq:gaussianint}.
\end{proof}

\section{Application to testing multiple hypotheses}
\label{sec:multiple}

In this appendix we apply the results of this paper
concerning multiple testing of a single hypothesis
to testing multiple hypotheses.
Namely, we will arrive at a generalization of the Bonferroni--Holm procedure \cite{Holm:1979}.
Fix a parameter
\begin{equation}\label{eq:parameter}
  r
  \le
  \frac{\ln K}{1-\ln K}
\end{equation}
(cf.\ \eqref{eq:min}); the Bonferroni--Holm case is $r:=-\infty$.

\begin{algorithm}[bt]
  \caption{Generalized Bonferroni--Holm procedure}
  \label{alg:GBH}
  \begin{algorithmic}
    \Require
      A significance level $\epsilon>0$ and parameter $r<-1$ (or, w.l.o.g., \eqref{eq:parameter}).
    \Require
      A sequence of p-values $p_1,\dots,p_K$ ordered as $p_{(1)}\le\dots\le p_{(K)}$.
    \For{$k=1,\dots,K$}
      \State $\texttt{reject}:=\texttt{true}$
      \State $I:=\{k\}$
      \For{$i=K,\dots,1,0$}
        \If{$\frac{r}{r+1}\left|I\right|^{1+1/r} M_{r,\left|I\right|}(p_{(I)}) > \epsilon$}
          \State $\texttt{reject}:=\texttt{false}$
        \EndIf
        \State $I:=I\cup\{i\}$
      \EndFor
      \If{$\texttt{reject}=\texttt{true}$}
        \State reject $H_k$
      \EndIf
    \EndFor
  \end{algorithmic}
\end{algorithm}

Suppose $p_k$ is a p-value for testing a composite null hypothesis $H_k$
(meaning that, for any $\epsilon\in(0,1)$, $\Prob(p_k\le\epsilon)\le\epsilon$ under $H_k$).
For $I\subseteq\{1,\dots,K\}$,
let $H_I$ be the hypothesis
\[
  H_I
  :=
  \left(
    \cap_{k\in I} H_k
  \right)
  \cap
  \left(
    \cap_{k\in\{1,\dots,K\}\setminus I} H^\complement_k
  \right),
\]
where $H^\complement_k$ is the complement of $H_k$.

Fix a significance level $\epsilon$.
Let us reject $H_I$ when
\[
  \frac{r}{r+1}\left|I\right|^{1+1/r}
  M_{r,\left|I\right|}(p_I)
  \le
  \epsilon,
\]
where $p_I$ is the vector of $p_k$ for $k\in I$;
by Proposition~\ref{prop:asymptotically-precise}, the probability of error will be at most $\epsilon$.
If we now reject $H_k$ when all $H_I$ with $I\supseteq\{k\}$ are rejected,
the family-wise error rate (FWER) will be at most $\epsilon$.
This gives the procedure described as Algorithm~\ref{alg:GBH},
in which $(k_1,\dots,k_K)$ stands for a fixed permutation of $\{1,\dots,K\}$
such that $p_{k_1}\le\dots\le p_{k_K}$.

\begin{algorithm}[bt]
  \caption{Generalized Bonferroni--Holm procedure for adjusting p-values}
  \label{alg:GBH-p}
  \begin{algorithmic}
    \Require
      A parameter $r<-1$ (or, w.l.o.g., \eqref{eq:parameter}).
    \Require
      A sequence of p-values $p_1,\dots,p_K$ ordered as $p_{(1)}\le\dots\le p_{(K)}$.
    \For{$k=1,\dots,K$}
      \State $p^*_k:=0$
      \State $I:=\{k\}$
      \For{$i=K,\dots,1,0$}
        \If{$\frac{r}{r+1}\left|I\right|^{1+1/r} M_{r,\left|I\right|}(p_{(I)}) > p^*_k$}
          \State $p^*_k := \frac{r}{r+1}\left|I\right|^{1+1/r} M_{r,\left|I\right|}(p_{(I)})$
        \EndIf
        \State $I:=I\cup\{i\}$
      \EndFor
    \EndFor
  \end{algorithmic}
\end{algorithm}

An alternative representation of the generalized Bonferroni--Holm procedure
given as Algorithm~\ref{alg:GBH}
is in terms of adjusting the p-values $p_1,\dots,p_K$ to new p-values $p^*_1,\dots,p^*_K$
that are valid in the sense of the FWER:
we are guaranteed to have $\Prob(\min_{k\in I}p^*_k\le\epsilon)\le\epsilon$ for all $\epsilon\in(0,1)$,
where $I$ is the set of the indices $k$ of the true hypotheses $H_k$.
The adjusted p-values can be defined as
\[
  p^*_k
  :=
  \max_{k\in I\subseteq\{1,\dots,K\}}
  \frac{r}{r+1}\left|I\right|^{1+1/r} M_{r,\left|I\right|}(p_I)
\]
and computed using Algorithm~\ref{alg:GBH-p}.

If we do not insist on controlling the FWER,
we can still use our ways of combining p-values instead of Bonferroni's
in more flexible procedures for testing multiple hypotheses,
such as those described in \cite{Goeman/Solari:2011}.
\end{document}